\documentclass[review,hidelinks,onefignum,onetabnum]{siamart220329}

\usepackage{hyperref}
\hypersetup{pdfauthor={Name}}
\hypersetup{colorlinks=true}

\usepackage{float}
\usepackage{amsfonts}
\usepackage{amsmath}
\usepackage{subcaption}
\usepackage{tikz}
\usetikzlibrary{trees,calc, positioning,decorations.markings,arrows,backgrounds,arrows.meta, bending,positioning,decorations.pathreplacing} 
\usepackage{graphicx}
\usepackage{xcolor}

\usepackage{algorithm}
\usepackage{algpseudocode}
\usepackage{makecell}

\usepackage{bm}         %
\usepackage{siunitx}    %
\usepackage{booktabs} %
\usepackage{multirow}   %
\usepackage{placeins} %
\usepackage{soul}

\usepackage[most]{tcolorbox}

\ifpdf
  \DeclareGraphicsExtensions{.eps,.pdf,.png,.jpg}
\else
  \DeclareGraphicsExtensions{.eps}
\fi

\usepackage{nameref}

\newtcbox{\highlight}[2]{enhanced, box align=base, nobeforeafter, colback=#1,
colframe=#2, size=small, left=0pt, right=0pt, boxsep=1.6pt, boxrule=1.2pt}

\newcommand{\githuburl}{https://github.com/Alexey-Voronin/HighOrderStokesMG}
\newcommand{\githuburltext}{https://github.com/Alexey-Voronin/HighOrderStokesMG}

\renewcommand{\vec}[1]{\bm{#1}}

\newcommand{\vecP}[1]{$\pmb{ \mathbb{P}}_{#1}$}
\newcommand{\Pk}{$\pmb{ \mathbb{P}}_k$}
\newcommand{\Pkmonedisc}{$\mathbb{P}_{k-1}^{disc}$}

\newcommand{\THPk}{$\pmb{\mathbb{P}}_{k}/\mathbb{P}_{k-1}$}
\newcommand{\THPkhat}{$\pmb{\mathbb{P}}_{\hat{k}}/\mathbb{P}_{\hat{k}-1}$}
\newcommand{\SVPk}{$\pmb{\mathbb{P}}_{k}/\mathbb{P}_{k-1}^{disc}$}
\newcommand{\SVPkhat}{$\pmb{\mathbb{P}}_{\hat{k}}/\mathbb{P}_{\hat{k}-1}^{disc}$}
\newcommand{\isoTHP}{$\pmb{\mathbb{P}}_1 \text{iso}\kern1pt\pmb{ \mathbb{P}}_2/ \mathbb{P}_1$}
\newcommand{\isoTHQ}{$\pmb{\mathbb{Q}}_1 \text{iso}\kern1pt\pmb{ \mathbb{Q}}_2/\mathbb{Q}_1$}
\newcounter{tempCounter}

\newcommand{\SVP}[1]{\setcounter{tempCounter}{#1}\addtocounter{tempCounter}{-1}$\pmb{\mathbb{P}}_{#1}/\mathbb{P}_{\thetempCounter}^{disc}$}
\newcommand{\THP}[1]{\setcounter{tempCounter}{#1}\addtocounter{tempCounter}{-1}$\pmb{\mathbb{P}}_{#1}/\mathbb{P}_{\thetempCounter}$}

\newcommand{\phmg}{$ph$\text{MG}}
\newcommand{\hmg}{$h$\text{MG}}
\newcommand{\afbf}{\text{FBF}}

\newcommand{\fhmg}{\afbf{}(\hmg{})}
\newcommand{\fphmg}{\afbf{}(\phmg{})}

\newcommand{\ASMVankaStarPC}{\text{ASMVankaStarPC}}
\newcommand{\ASMStarPC}{\text{ASMStarPC}}
\newcommand{\PCPatch}{\text{PCPatch}}

\newcommand{\VphmgDirect}{$ph$\text{MG(direct)}}
\newcommand{\VphmgGradual}{$ph$\text{MG(gradual)}}

\newcommand{\Vphmg}[4]{$V^{(n_V=#3,
\ifnum#1=#2
    \nu_p=\nu_h=#1
\else
    \nu_p=#1,\nu_h=#2
\fi
)}_{ph\text{MG(#4)}}$
}

\newcommand{\relSetupFull}{\textbf{Rel. T(setup)}}
\newcommand{\relTotalFull}{\textbf{Rel. T(total)}}
\newcommand{\relSolveFull}{\textbf{Rel. T(solve)}}
\newcommand{\setupOverTotalFull}{\textbf{T(setup)/T(total)}}

\newcommand{\relSetup}{\textbf{R(setup)}}
\newcommand{\relTotal}{\textbf{R(total)}}
\newcommand{\relSolve}{\textbf{R(solve)}}
\newcommand{\setupOverTotal}{\textbf{T(setup)}}

\definecolor{tab:green}{RGB}{0,176,80}
\definecolor{tab:red}{RGB}{220,0,0} 

\Crefname{algocf}{Algorithm}{Algorithms}
\newtheorem{remark}[theorem]{Remark}

\newcommand{\FULLTITLE}{Monolithic Multigrid Preconditioners for High-Order Discretizations of Stokes
 Equations}
\newcommand{\SHORTTITLE}{Monolithic MG for High-Order Stokes}

\title{\FULLTITLE\thanks{%
}}
\headers{\SHORTTITLE}{A. Voronin, G. Harper, S. MacLachlan, L.N. Olson, and R.S. Tuminaro}

\author{Alexey Voronin\thanks{Department of Computer Science, University of
  Illinois at Urbana-Champaign, Urbana IL, 61801 USA, and Sandia National
  Laboratories, Albuquerque, NM 87123 USA. (\email{abvoron@sandia.gov}, \url{http://alexey-voronin.github.io/}).}
 \and Graham Harper\thanks{Sandia National Laboratories, Livermore, CA 94551
 USA (\email{gbharpe@sandia.gov}, \email{rstumin@sandia.gov}),
 (\url{https://www.sandia.gov/ccr/staff/graham-bennett-harper/}, \url{https://www.sandia.gov/ccr/staff/raymond-s-tuminaro/}).
   }
 \and Scott MacLachlan\thanks{Department of
  Mathematics and Statistics, Memorial University of Newfoundland, St.~John's NL, Canada
   (\email{smaclachlan@mun.ca}, \url{https://www.math.mun.ca/\string~smaclachlan/}).
   }
 \and Luke N. Olson\thanks{Department of Computer Science, University of
   Illinois at Urbana-Champaign, Urbana IL
   (\email{lukeo@illinois.edu}, \url{http://lukeo.cs.illinois.edu/}).}
 \and Raymond S. Tuminaro\footnotemark[3]
 }

\usepackage{amsopn}

\ifpdf
\hypersetup{
  pdftitle={\FULLTITLE},
  pdfauthor={A. Voronin, S. MacLachlan, L. N. Olson, R. Tuminaro, and G. Harper}
}
\fi

\begin{document}

\maketitle

\begin{abstract}
This work introduces and assesses the efficiency of a monolithic \phmg{} multigrid framework designed for high-order discretizations of stationary Stokes systems using Taylor-Hood and Scott-Vogelius elements. 
The proposed approach integrates coarsening in both approximation order ($p$) and mesh resolution ($h$), to address the computational and memory efficiency challenges that are often encountered in conventional high-order numerical simulations. 
Our numerical results reveal that \phmg{} offers significant improvements over traditional spatial-coarsening-only multigrid (\hmg{}) techniques for problems discretized with Taylor-Hood elements across a variety of problem sizes and discretization orders. In particular, the \phmg{} method exhibits superior performance in reducing setup and solve times, particularly when dealing with higher discretization orders and unstructured problem domains. 
For Scott-Vogelius discretizations, while monolithic \phmg{} delivers low
    iteration counts and competitive solve phase timings, it exhibits a
    discernibly slower setup phase when compared to a multilevel
    (non-monolithic) full-block-factorization (\afbf{}) preconditioner where
    \phmg{} is employed only for the velocity unknowns. This is primarily due to
    the setup costs of the larger mixed-field relaxation patches with monolithic
    \phmg{} versus the patch setup costs with a single unknown type for \afbf{}.
\end{abstract}



\section{Introduction}

High-order discretization methods are increasingly favored in engineering for
their numerical accuracy~\cite{demkowicz1989toward,rachowicz1989toward}, reduced
numerical dissipation, and improved efficiency in parallel computing
environments~\cite{fischer2020scalability,kolev2021efficient}. Despite these
advantages, a notable drawback is the poor conditioning of the resulting linear and linearized systems of equations,
leading to reduced convergence rates for many iterative linear solvers. To address this,
our study aims to identify effective, efficient, and scalable preconditioners
for stationary Stokes systems, adaptable across different discretization orders
of Taylor-Hood (\THPk{}) and Scott-Vogelius discretizations (\SVPk{}).

We explore spatial-coarsening $h$-multigrid (\hmg{}) and
$p$\kern1pt-multigrid methods that coarsen in approximation order for Stokes
systems. These methods draw on insights from the discretization of the Poisson equation, where suitably defined
\hmg{} is known to provide optimal convergence rates irrespective of discretization
order or mesh size. However, the memory footprint of \hmg{} at higher
discretization orders limits the size of problems that can be pracitcally solved.
An effective strategy that helps mitigate this is coarsening in order ($p$),
followed by spatial coarsening, a technique that significantly reduces memory
requirements on coarser grids. This approach, denoted as \phmg{}, conserves
memory and minimizes computational and communication demands at higher
discretization orders~\cite{tielen2020p,fischer2015scaling,voronin2022low}.
Nonetheless, \phmg{} can sometimes face challenges in achieving robust convergence,
particularly with aggressive coarsening strategies~\cite{thompson2023local}.
This paper extends \phmg{} to Stokes systems, identifying the most effective solver strategies for various discretizations, orders, and problem sizes.

Other strategies for reducing memory footprint include utilizing hybridizable
discretizations with static condensation~\cite{he2021local}, low-order-refined
preconditioning~\cite{orszag1979spectral,canuto2010finite,casarin1997quasi,pazner2023end,pazner2024matrix}
and $p$\kern1pt-multigrid~\cite{botti2022p}.
These approaches can require extensive efforts in code development and
discretization formulation, including the careful tuning of stability parameters
to ensure convergence and numerical stability.
Geometric multigrid~\cite{saye2020fast} and matrix-free geometric
multigrid methods for Stokes~\cite{kohl2022textbook,jodlbauer2022matrix}
also reduce the memory footprint,
but have not been demonstrated in literature for the large polynomial degrees we consider.

A critical factor in the robustness of both \hmg{} and \phmg{} methods is the selection of relaxation schemes. Our focus in this study is the additive variant of the Schwarz relaxation method, also known as overlapping patch relaxation or (additive) Vanka relaxation. The effectiveness of these relaxation methods relies on a decomposition of the domain into subdomains in a way that aligns with the structure of the discretization~\cite{farrell2021pcpatch,arnold1997preconditioning,farrell2021reynolds,arnold2000multigrid}. Given a sound decomposition, these relaxation methods achieve robust convergence in both scalar and coupled linear systems~\cite{schoberl2003schwarz,farrell2021local} and are more scalable on parallel architectures than their multiplicative counterpart~\cite{farrell2021local}. 
Multiple domain decompositions are sometimes possible for a given discretization, each with its unique structure, size, and convergence properties. We discuss several of these decompositions and highlight their performance in high-order monolithic multigrid solvers.

Additionally, using \phmg{} over \hmg{} allows for the integration of algebraic multigrid (AMG) techniques at the coarsest discretization level. 
Considering the robustness challenges associated with directly applying AMG to high-order discretizations,
employing the \phmg{} framework still permits the utilization of AMG on the
coarsest-in-order discretization within the \phmg{} solver. This should prove particularly
advantageous for complex domains and unstructured
meshes~\cite{voronin2023monolithic,notay2017algebraic,bacq2022new}, although we do not investigate this use of AMG here.

Recent years have seen significant interest in the development of multilevel solvers for high-order discretizations of the Stokes equations for both Taylor-Hood and Scott-Vogelius elements~\cite{jodlbauer2022matrix, botti2022p, voronin2023monolithic, saye2020fast, farrell2021reynolds, kohl2022textbook, Rafiei_2024, bacq2022new}.  In this work, we examine the impacts of how the multilevel hierarchy is constructed on the performance of the resulting solver.  We demonstrate that relatively rapid coarsening in order leads to effective preconditioners for both discretizations, and compare the performance of
\phmg{} and \hmg{} preconditioners, highlighting key factors that influence
their effectiveness.  For cases where high discretization orders are required, we offer insight into the computational and memory efficiency challenges that arise, and how they can be addressed.

We note that the focus of this work is on the efficiency of multilevel solvers for high-order Stokes
discretizations, but that the choice between high-order and low-order
methods should be informed by analysis that weighs the superior approximation
quality of high-order schemes against the extra computational costs per degree
of freedom. This paper does not offer specific recommendations for choosing the
discretization order or type, as these decisions are problem-dependent. For
further guidance on this matter, refer
to~\cite{farrell2021reynolds,case2011connection,betteridge2021code,linke2011convergence}.

The remainder of this paper is structured as follows. In~\cref{sec:disc}, we introduce the Stokes equations and the accompanying discretizations. \Cref{sec:precodnitioners} describes the \hmg{} and \phmg{} frameworks for monolithic and block preconditioners of Stokes systems, detailing the patch relaxation approaches used for each preconditioner type. \Cref{sec:num_results} presents numerical results, demonstrating the robustness of the preconditioner for a set of test problems and comparing the computational cost across them. \Cref{sec:lessons_learned} summarizes the insights gained from the numerical experiments. Finally, \cref{sec:conclusion} offers conclusions, and \cref{sec:future_directions} outlines future research directions.

\section{Discretization of the Stokes Equations}\label{sec:disc}

This paper considers the steady-state Stokes problem on a bounded Lipschitz domain $\Omega \subset \mathbb{R}^d$, where $d=2$ or $3$. The boundary of this domain is denoted as $\partial\Omega=\Gamma_{\text{N}} \cup \Gamma_{\text{D}}$, and the problem is formulated as follows:
\begin{subequations}\label{eq:stokes-eq}
  \begin{alignat}{2}
    -\nabla^2 \vec{u} +\nabla p &=\vec{f}              &\quad& \text{in $\Omega$} \label{eq:stokes-eq1} \\
    -\nabla \cdot \vec{u}       &= 0                   &\quad& \text{in $\Omega$} \label{eq:stokes-eq2} \\
                        \vec{u} &= \vec{g}_{\text{D}}  &\quad& \text{on $\Gamma_{\text{D}}$} \label{eq:stokes-eq3} \\
        \frac{\partial\vec{u}}{\partial \vec{n}} - \vec{n} p
                                &= \vec{g}_{\text{N}}  &\quad& \text{on $\Gamma_{\text{N}}$}. \label{eq:stokes-eq4}
\end{alignat}
\end{subequations}
In these equations, $\vec{u}$ represents the velocity of a viscous fluid, $p$ is the pressure, $\vec{f}$ is a forcing term, and $\vec{n}$ is the outward pointing normal to $\partial\Omega$. The terms $\vec{g}_{\text{D}}$ and $\vec{g}_{\text{N}}$ denote the given boundary data, while $\Gamma_{\text{D}}$ and $\Gamma_{\text{N}}$ refer to disjoint Dirichlet and Neumann segments of the boundary, respectively. In a variational formulation, the natural function space for the velocity is $\vec{\mathcal{H}}_{\vec{g}_{\text{D}}}^1(\Omega)$, a subset of $\vec{\mathcal{H}}^1(\Omega)$ that satisfies the essential boundary condition in~\eqref{eq:stokes-eq3}.
If the velocity is specified everywhere along the boundary, which necessitates $\int_{\partial\Omega} \vec{g}_{\text{D}}\cdot\vec{n}=0$, the appropriate function space for the pressure is the set of zero-mean functions in $L^2(\Omega)$, denoted as $\mathcal{Q}=L^2_0(\Omega)$. This choice helps obtain a unique pressure solution, not merely up to an additive constant. In scenarios where the measure of $\Gamma_{\text{N}}$ is non-zero, this additional constraint becomes redundant, and the pressure space is then assumed to be $\mathcal{Q}=L^2(\Omega)$~\cite{elman2014finite}.

We denote the
finite-dimensional spaces associated with a mixed finite element discretization %
by $(\vec{\mathcal{V}}_h^{\text{D}}, \mathcal{Q}_h)\subset (\vec{\mathcal{H}}_{\vec{g}_{\text{D}}}^1(\Omega), \mathcal{Q})$. Here, $\vec{\mathcal{V}}_h^{\text{D}}$ strongly satisfies the Dirichlet boundary conditions on $\Gamma_{\text{D}}$ as outlined in~\cref{eq:stokes-eq3} up to interpolation error. The weak formulation derived from~\cref{eq:stokes-eq} seeks $\vec{u} \in \vec{\mathcal{V}}_h^{\text{D}}$ and $p \in \mathcal{Q}_h$, satisfying the following system for all $\vec{v} \in \vec{\mathcal{V}}_h^0$ and $q \in \mathcal{Q}_h$, where $\vec{\mathcal{V}}_h^0$ is akin to $\vec{\mathcal{V}}_h^{\text{D}}$, but with zero Dirichlet boundary conditions on $\Gamma_{\text{D}}$:
\begin{subequations}\label{eq:stokes-disc}
\begin{alignat}{2}
& a(\vec{u}, \vec{v})+b(p, \vec{v}) &&= F(\vec{v}) \label{eq:weakform1}\\
& b(q, \vec{u}) &&= 0 \label{eq:weakform2},
\end{alignat}
\end{subequations}
where $a(\cdot,\cdot)$ and $b(\cdot,\cdot)$ represent bilinear forms, and $F(\cdot)$ is a linear form defined as
\begin{equation*}
a(\vec{u}, \vec{v}) = \int_{\Omega} \nabla \vec{u} : \nabla \vec{v},\quad
b(p, \vec{v}) = - \int_{\Omega} p \nabla \cdot \vec{v},\quad\text{and}\quad
F(\vec{v}) = \int_{\Omega} \vec{f} \cdot \vec{v}+\int_{\Gamma_{\text{N}}} \vec{g}_{\text{N}} \cdot \vec{v}.
\end{equation*}
Henceforth, the notation $\vec{u}$ and $p$ will also represent the discrete velocity and pressure unknowns in the finite-element discretization of~\cref{eq:stokes-disc}. Given an inf-sup stable choice of finite-dimensional spaces $(\vec{\mathcal{V}}_h^{\text{D}}, \mathcal{Q}_h)$, a saddle-point system emerges, as shown in equation~\eqref{eq:saddle}. The matrices $A$ and $B$ in this system represent the discrete vector Laplacian on $\mathcal{V}_h^{\text{D}}$ and the negative of the discrete divergence operator mapping $\mathcal{V}_h^{\text{D}}$ into $\mathcal{Q}_h$, respectively.
\begin{equation}\label{eq:saddle}
    Kx
    =
    \begin{bmatrix}
        A  & B^T \\
        B  &  0
    \end{bmatrix}
    \begin{bmatrix}
        \vec{u}\\
        p
    \end{bmatrix}
    =
    \begin{bmatrix}
        \vec{f}\\
        0
    \end{bmatrix}
    =b,
\end{equation}
On $d$ dimensional domains where $\vec{\mathcal{V}}_h^{\text{D}}$ is a vector-valued version of a scalar function space (such as the vector Lagrange elements considered here), the dimensions of matrices $A$ and $B$ can be expressed as $A \in \mathbb{R}^{(d\times n_{u})\times (d\times n_u)}$ and $B \in \mathbb{R}^{n_{p}\times (d\times n_{u})}$. Here, $n_u$ is the number of velocity nodes for each component of the vector $\vec{u}$, and the term $n_p$ represents the number of pressure nodes.

\subsection{Finite Element Discretizations and Meshes}
\begin{figure}[!ht]
  \centering
  \includegraphics[width=\textwidth]{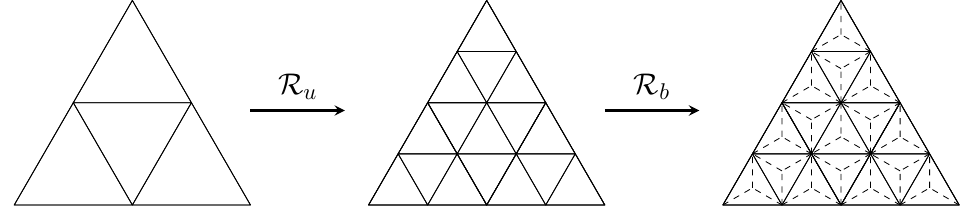}
  \caption{Sequential Mesh Refinement Processes: The initial image (left) illustrates the triangulation of a triangular domain. This is followed by a quadrisection refinement ($\mathcal{R}_u$), resulting in the intermediate mesh (middle). The final image (right) shows the result of a barycentric mesh refinement ($\mathcal{R}_b$), employing the Alfeld split, further refining the intermediate mesh.}
\label{fig:meshes}
\end{figure}

We focus on two distinct but related discretizations employed to define the finite-dimensional spaces $(\vec{\mathcal{V}}_h^{\text{D}}, \mathcal{Q}_h)$ on meshes composed of $d$-dimensional simplices. The first discretization we consider is the well-known Taylor-Hood (TH) discretization, \THPk{}~\cite{elman2014finite}. This discretization utilizes continuous piecewise polynomials of degree $k$ for the velocity space $\vec{\mathcal{V}}_h^{\text{D}}$ and continuous piecewise polynomials of degree $k-1$ for the pressure space $\mathcal{Q}_h$. Notably, the TH element pair is inf-sup stable for $k \geq 2$ on any triangular (tetrahedral) mesh $\Omega_h$ of the domain $\Omega$ as depicted in~\cref{fig:meshes}. 
The Scott-Vogelius (SV) discretization, \SVPk{}, employs continuous piecewise polynomials of degree $k$ for the velocity and discontinuous piecewise polynomials of degree $k-1$ for the pressure. The inf-sup stability of the SV discretization depends on the polynomial degree, $k$, of the function spaces and the characteristics of the mesh utilized in assembling the system. As demonstrated in~\cite{qin1994convergence,zhang2005new}, the \SVPk{} discretization achieves stability for $k \geq d$ on any mesh resulting from a single step of barycentric refinement of any triangular (tetrahedral) mesh. An example of such mesh refinements is depicted in middle and right meshes in~\cref{fig:meshes}, where $\mathcal{R}_b$ symbolizes the barycentric refinement of $\mathcal{R}_u(\Omega_h)$. Further, studies in~\cite{zhang2008p1,zhang2011quadratic} have shown that for meshes derived from a Powell-Sabin split of a triangular or tetrahedral mesh, the polynomial order required for inf-sup stability in the SV discretization is reduced to $k \geq d-1$. However, such meshes are not within the scope of this work.

A notable advantage of the SV discretization over the TH discretization is its ability to ensure that the divergence of any velocity in \Pk{} is precisely represented in the pressure space, \Pkmonedisc{}. This leads to velocities that are point-wise divergence-free when the divergence constraint is enforced weakly.  Notably, TH elements do not have this property, although the addition of a grad-div stabilization term can enhance point-wise mass conservation, as discussed in~\cite{john2017divergence,case2011connection,linke2011convergence}. The point-wise divergence-free property of the SV discretization comes at the increased cost of a much larger number of unknowns in the pressure field relative to the same size mesh and discretization order for TH. The discontinuous pressure field also introduces new challenges for designing efficient multilevel iterative methods, most of which have been designed with continuous bases in mind.

\section{Preconditioners}\label{sec:precodnitioners}

This section introduces several multilevel algorithms for preconditioning Stokes systems, $Kx=b$, as defined in~\cref{eq:saddle}, for Taylor-Hood and  Scott-Vogelius discretizations. 
The motivation for the development of these preconditioners stems from the observation that \hmg{} methods can achieve robust convergence rates across any discretization and problem size on structured meshes~\cite{Rafiei_2024}. However, the computational cost associated with these methods becomes prohibitively high at high discretization orders.  This is true already for geometric \hmg{} methods, but becomes even worse when algebraic multigrid is applied directly to higher-order discretizations, primarily due to the growth of the stencil size on coarse grids.
To address this challenge, we investigate employing lower-order discretizations on coarser grids to achieve similar convergence rates while significantly reducing computational expenses. We aim to identify when \hmg{} should be preferred over \phmg{}, especially considering the setup cost and time to solution. Since the cost of both types of solvers is largely dominated by relaxation, we examine approaches that minimize relaxation costs across all multigrid levels.

In~\cref{sec:hmg}, we describe the assembly of monolithic hierarchical multigrid
(\hmg{}) methods for Taylor-Hood discretizations. In~\cref{sec:phmg}, we
introduce a ``coarsening in order'' multigrid approach, \phmg{}, transitioning
from high- to low-order discretizations to improve solver efficiency.
In~\cref{sec:afbf}, we discuss the full-block factorization preconditioner
(\afbf{}) for Scott-Vogelius discretizations, where standard \hmg{} does not yield an effective solver.
The section concludes with~\cref{sec:rlx}, describing patch relaxation
approaches for each preconditioner.

\subsection{\protect\hmg{}}\label{sec:hmg}

\begin{figure}[!htp]
\centering
\includegraphics[width=0.55\linewidth]{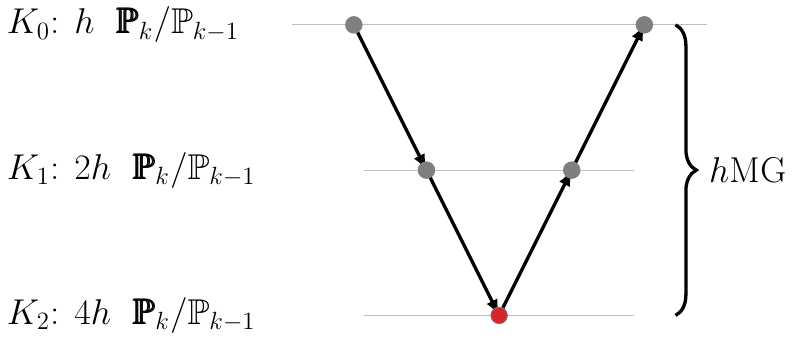}
\caption{Illustration of the \protect\hmg{} cycle, showcasing the rediscretization process on coarser meshes with the same discretization order. The relaxation stages are represented by gray dots at each level of the V-cycle. Black lines signify the grid-transfer points, while the red dot at the base of the V-cycle indicates an exact coarse-grid solve.}
\label{fig:hmg_cycle}
\end{figure}

This study employs geometric multigrid (GMG) methods with standard coarsening, characterized by doubling the grid size $h$ at each hierarchy level. %
Central to our approach is the two-level error-propagation operator for a multigrid cycle, described as:
\begin{equation}
G_{\ell} = (I - \omega_{\ell} M_{\ell}^{-1} K_{\ell})^{\nu_2} (I - P_{\ell}
K_{\ell+1}^{-1} R_{\ell} K_{\ell}) ( I -\omega_{\ell} M_{\ell}^{-1}
K_{\ell})^{\nu_1} \text{ for } \ell>0, \label{eq:TG-operator}
\end{equation}
Here, $\ell$ represents the current level in the multigrid hierarchy. $K_{\ell}$ is the system matrix at level $\ell$, and $M_{\ell}$ is the corresponding relaxation operator, with $\omega_{\ell}$ as the damping factor. The terms $\nu_1$ and $\nu_2$ represent the numbers of pre-relaxation and post-relaxation sweeps, respectively. $P_{\ell}$ and $R_{\ell}$ are the interpolation and restriction operators for transitions between hierarchy levels. In our approach, we use finite-element interpolation for $P_{\ell}$, and the restriction operators $R_{\ell}$ are defined as the transpose of the interpolation operators. In GMG methods, $K_{\ell+1}$ is typically a rediscretization of the original problem on a coarser mesh. However, an equivalent coarse grid can be computed via the Galerkin product $K_{{\ell}+1}=R_{\ell} K_{\ell} P_{\ell}$, at least when no stabilization terms are present (as in the cases considered here).
The $G_0$ in Equation~\cref{eq:TG-operator} represents a two-level solver for approximating solutions to the $K_0$ system. Equation~\cref{eq:TG-operator} can be recursively extended to three-level \hmg{} which we depict in~\cref{fig:hmg_cycle}.

The mechanism of assembly of interpolation operators for each scalar field is equivalent to $L_2$ projection of $u_k \in \mathcal{V}_{k} $ onto $\mathcal{V}_{k+1}$. This operation is defined as finding $u_{k+1}\in\mathcal{V}_{k+1}$ such that  $\langle u_{k+1}, v_{k+1}\rangle=\langle u_{k}, v_{k+1}\rangle$ for all $v_{k+1}\in\mathcal{V}_{k+1}$. 
As $u_{k+1}$ and $u_{k}$ can be expressed in terms of basis functions, this leads to solving a matrix equation involving a standard $\mathcal{V}_{k+1}$ mass matrix, where the right-hand side is defined in terms of basis functions of the two spaces. In the context of $hMG$, this can be accomplished without the solution of any linear systems for nested meshes. The restriction operator is obtained as a transpose of the resulting interpolation operator.

\subsection{\protect\phmg{}}\label{sec:phmg}

This section introduces the defect-correction approach, $G_\text{dc}$, as a
cornerstone of our $p$\kern1pt-coarsening strategy within a geometric multigrid
framework. We employ two distinct Stokes operators for this purpose: the
high-order discretizations $K_0$ (\THPk{} or \SVPk{}, for $k \geq 3$) and the
simpler, lower-order discretization $K_1$ (\THP{2}). This double discretization
strategy will be shown to improve overall computational efficiency, by substituting high-order
computations with more economical lower-order approximations on coarse grids without substantially 
compromising the overall convergence rates of our multilevel solvers. Key to
integrating these discretizations are the transfer operators $P_0$ and
$R_0$, which map functions between the high and low-order bases, thereby
maintaining the effectiveness of the solvers across different resolutions. These
grid-transfer operators are assembled similarly to~\hmg{} grid-transfer
operators, as described in~\cref{sec:hmg}.

\begin{figure}[!ht]
\centering
\begin{subfigure}{0.49\textwidth}
 \includegraphics[width=\textwidth]{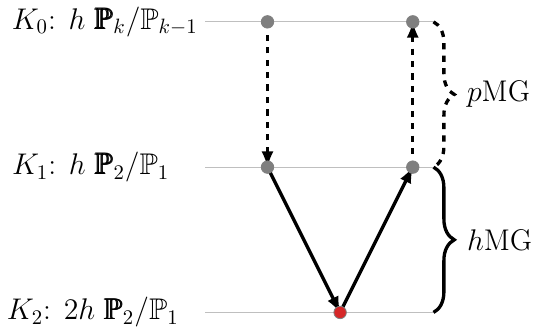}
 \caption{Direct \protect\phmg{}}\label{fig:phmg_direct}
\end{subfigure}
\hfill
\begin{subfigure}{0.49\textwidth}
 \includegraphics[width=\textwidth]{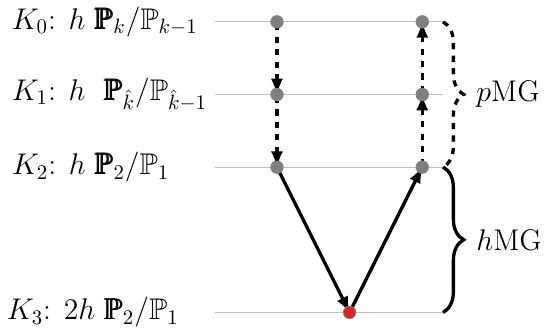}
 \caption{Gradual \protect\phmg{}}\label{fig:phmg_gradual}
\end{subfigure}
\caption{Illustration of \protect\phmg{}  V-cycles applied to \protect\THPk{} discretizations,  deviating from~\protect\cref{fig:hmg_cycle} only in the modal restriction and interpolation represented by dashed lines. These cycles are also compatible with Scott-Vogelius discretizations (\protect\SVPk{}) on $\ell=0$. Figure (a) depicts a Direct \protect\phmg{} cycle where the \protect\THPk{} discretization directly transitions to \protect\THP{2}. Figure (b) shows a Gradual \protect\phmg{} cycle, where \protect\THPk{} discretization first steps down to an intermediate \protect\THPkhat{} with $k > \hat{k} > 2$, before reducing the approximation order to \protect\THP{2}.}
\label{fig:phmg_cycles}
\end{figure}

The error-propagation operator $G_\text{dc}$ %
for the defect-correction approach applied to $K_0$ is %
\begin{equation}\label{eq:DC-TG-operator}
G_\text{dc} =\left(I - \omega_0 M_0^{-1} K_0\right)^{{\nu}_2}
      \left(I - {P}_0(I-G_1) K_1^{-1} {R}_0 K_0\right)
      \left(I - \omega_0 M_0^{-1} K_0\right)^{{\nu}_1}.
\end{equation}
When $G_1 = 0$ (corresponding to an exact coarse solve using $K_1$), this corresponds to the two-level solver~\eqref{eq:TG-operator}
using a different coarse discretization and different grid transfer operators.
Here, ${\nu}_1$ and ${\nu}_2$ represent the counts of pre- and post-relaxation sweeps on $K_0$, where $M_0^{-1}$ denotes the relaxation scheme and $\omega_0$ denotes a relaxation weight. \Cref{eq:DC-TG-operator} consists of relaxation on $K_0$, a subsequent coarse/low-order correction ($G_1$), and then final post-relaxation on $K_0$. In this work, the low-order correction operator, $G_{1}$, 
will primarily be \hmg{}, as defined in~\cref{eq:TG-operator}. However, it can be replaced by a recursive application of~\cref{eq:DC-TG-operator}. This recursive application is depicted in~\cref{fig:phmg_cycles}, where~\cref{fig:phmg_direct,fig:phmg_gradual} illustrate the ``direct'' and ``gradual'' \phmg{} cycles, respectively. These figures show how \THPk{} discretization is coarsened, either transitioning directly to \THP{2} or stepping down through an intermediate \THPkhat{} stage. While either an intermediate \THPkhat{} or an intermediate \SVPkhat{} could be used when $K_0$ corresponds to the \SVPk{} discretization, we only gradually coarsen via \THPkhat{} as we have observed this yields similar convergence rates and cheaper relaxation costs. 

This approach builds on~\cite{voronin2022low,voronin2023monolithic} where we utilized the low-order~\isoTHP{} finite-element discretization to develop robust AMG preconditioners for more complex Stokes equations discretizations, such as \THP{2} and \SVP{2}. %
Although the \isoTHP{} discretization could be used here as well, an efficient parallel implementation of solvers for this discretization are not currently available in our solver software, so are not considered here.

\subsection{Approximate Full-Block Factorization (\protect\afbf{})}\label{sec:afbf}

Here, we present a full-block factorization preconditioner, as
discussed in~\cite{benzi2005numerical,elman2008taxonomy} for a general setting and in~\cite{farrell2021reynolds} 
for the Scott-Vogelius discretization of the Navier-Stokes equations. The \afbf{} preconditioner is 
\begin{equation}
P^{-1} =
\left( \begin{array}{cc}
I & -\tilde{A}^{-1} B^\top \\
0 & I \\
\end{array} \right)
\left( \begin{array}{cc}
\tilde{A}^{-1} & 0 \\
0 & \tilde{S}^{-1} \\
\end{array} \right)
\left( \begin{array}{cc}
I & 0 \\
-B\tilde{A}^{-1} & I \\
\end{array} \right),
\label{eq:fbf}
\end{equation}
Here, $\tilde{S}^{-1}$ is the inverse of the block diagonal pressure mass matrix, which is 
spectrally equivalent to the inverse Schur complement~\cite{verfurth1984combined}. The term $\tilde{A}^{-1}$ denotes the application of a \phmg{} V-cycle to $A$. %

Our version of the \afbf{} preconditioner differs from the ``ALFI'' approach in~\cite{farrell2021reynolds}
which 
targets the stationary
incompressible Navier-Stokes equations at high Reynolds numbers. The ALFI method
employs an augmented-Lagrangian technique coupled with relatively large relaxation patches
to achieve solver robustness for large Reynolds number. %
In contrast, solvers for viscous Stokes flow do not require
the augmented-Lagrangian technique, achieving robust convergence with more
compact patches. Additionally, ALFI utilizes a hierarchy of
non-nested, barycentrically refined meshes requiring customized interpolation
operators. Our approach avoids this by employing standard mesh refinement for
coarser grids and reserving barycentric refinement for the finest mesh only, as
illustrated in~\cref{fig:meshes}. This streamlined approach yields effective
preconditioning for the Stokes equations, but avoids the computational intricacies needed for the high Reynolds number setting for which ALFI was proposed.

\subsection{Relaxation}\label{sec:rlx}

Classical point relaxation methods such as Jacobi and Gauss-Seidel often fail to define robust relaxation methods for coupled PDE systems such as Stokes. 
To overcome this challenge, subspace-correction methods~\cite{xu1992iterative} have grown in popularity as relaxation schemes for these systems, due to their robustness, relative ease of implementation, and parallelizability.  
These methods leverage the decomposition of the solution space, $V_h = \sum_i V^{(i)}_h$ (where this need not be a direct sum) into subspaces $V^{(i)}_h$. 
For each subspace, $V^{(i)}_h$, an injection operator, $I^{(i)}$ is defined and used to create a patch submatrix, $K^{(i)}=I^{(i)} K  ({I^{(i)}})^T$. 
Together with a partition of unity operator, $W_i$, which addresses overlap in the space decomposition, these patches form the basis for a single iteration of the additive Schwarz method (ASM), written as a preconditioner as
\begin{equation}\label{eq:ASM}
M^{-1} = \sum_i {I^{(i)}}^T W_i \left(K^{(i)}\right)^{-1}  I^{(i)}.
\end{equation}
For a single fixed-point ASM iteration, the relaxation operator is $I-\omega M^{-1} K$, with $\omega$ being a damping parameter.

The ASM approach in~\cref{eq:ASM} is essentially a block Jacobi relaxation scheme, which retains Jacobi's suitability for parallel computation and efficacy in multigrid relaxation. 
Alternatively, one could use block Gauss-Seidel (multiplicative Schwarz) methods that, while less parallelizable, have demonstrated strong convergence rates for saddle-point systems~\cite{wang2013multigrid,farrell2021local}.
The effectiveness of both additive and multiplicative Schwarz methods relies heavily on the choice of space decomposition. This aspect has been thoroughly explored in previous studies~\cite{arnold1997preconditioning,arnold2000multigrid,farrell2021pcpatch,farrell2021reynolds}.

\begin{figure}[!ht]
\centering
\begin{subfigure}{0.44\textwidth}
\includegraphics[width=\textwidth]{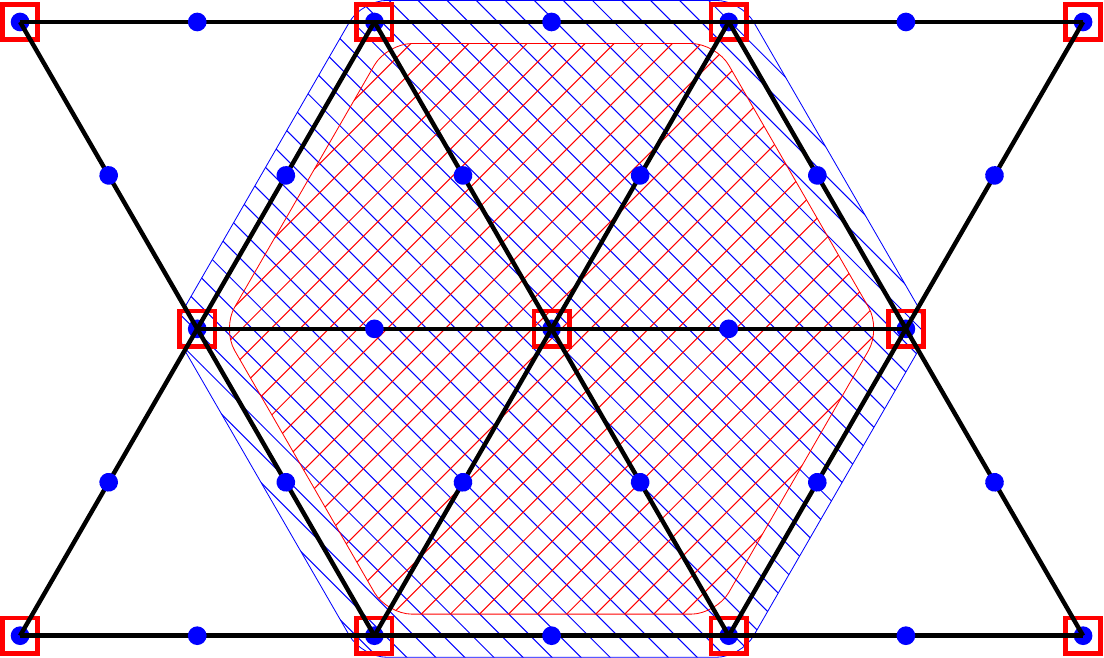}
\caption{Taylor-Hood patch}
\label{fig:th_patch}
\end{subfigure}
\hfill
\begin{subfigure}{0.53\textwidth}
\includegraphics[width=\textwidth]{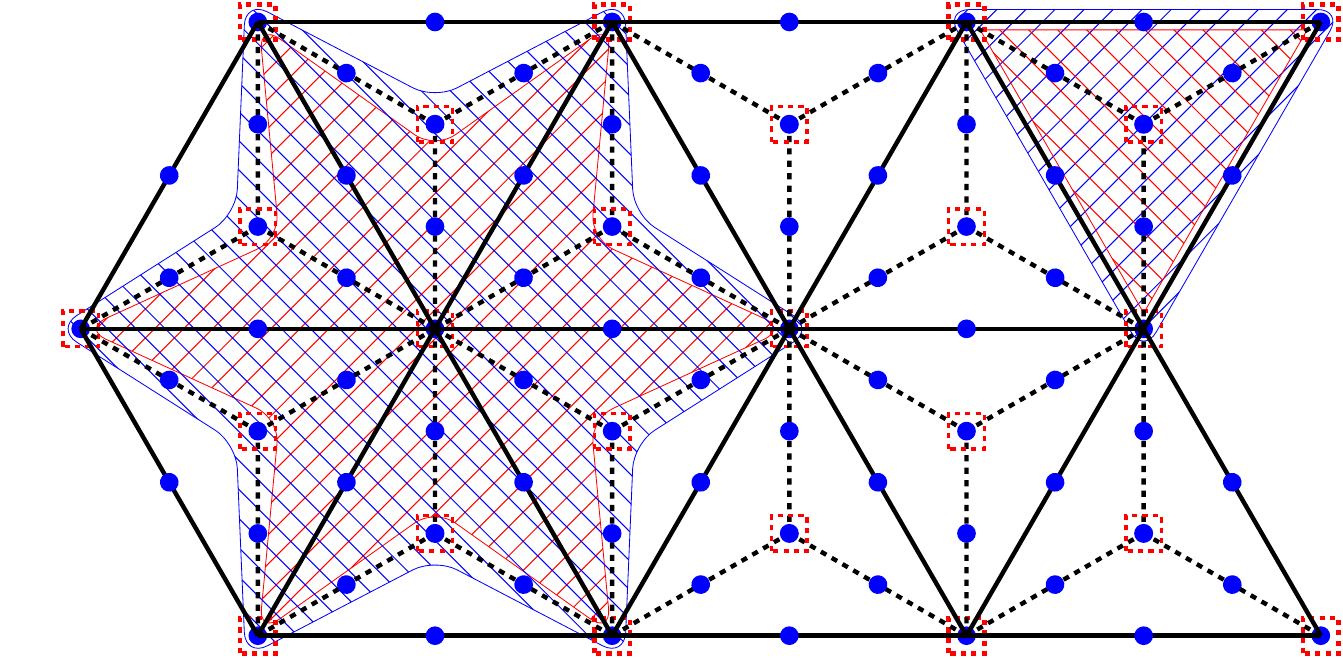}
\caption{Scott-Vogelius patch}
\label{fig:sv_patch}
\end{subfigure}     
\caption{%
	Blue dots (\protect\tikz[baseline=-0.5ex]\protect\draw[blue, fill=blue] (0, 0) circle (2pt);) represent velocity DoFs in $\mathbb{P}_2$.
	Red squares (\protect\tikz[baseline=-0.7ex]\protect\draw[line width=0.05mm, solid, thick, red] (0,0) +(-0.1,-0.1) rectangle ++(0.1,0.1);) correspond to pressure DoFs in $\mathbb{P}_1$. 
	Red framed squares with dashed lines (\protect\tikz[baseline=-0.7ex]\protect\draw[line width=0.08mm, dash pattern=on 1.pt off 1.pt, thick, red] (0,0) +(-0.1,-0.1) rectangle ++(0.1,0.1);) denote pressure DoFs in $\mathbb{P}_1^{disc}$.  
	Areas shaded red with northwest lines highlight the topological vertex star operator. Blue areas with northeast lines indicate the topological closure of the vertex star operator.
	\afbf{} employs \ASMStarPC{} patches, which include the star operation at each mesh vertex encompassing the enclosed velocity DoFs.
	For \phmg{}, \ASMVankaStarPC{} patches are used, which comprise the closure of the star operation for velocity DoFs and the star operation for pressure. In Scott-Vogelius discretizations, the star operation encompasses the discontinuous pressure DoFs associated with elements adjacent to the center vertex. The monolithic multigrid patches for Scott-Vogelius discretization are of two types: barycentric vertex patches (right) and other vertex patches (left).}
\label{fig:patches}
\end{figure}

When using monolithic multigrid (and not the FBF preconditioner), we adopt a coupled relaxation approach for the \hmg{} and \phmg{} TH and SV
preconditioners, where both velocity and pressure DoFs are relaxed
simultaneously. At each mesh vertex, we form patches using a vertex-star topological operation, which selects all degrees-of-freedom (DoFs) on topological entities (edges, faces, and elements) adjacent to the given vertex~\cite{farrell2021pcpatch}.
Specifically, for a given mesh vertex, the pressure DoFs in the resulting patch are those contained in the star of the vertex, while the velocity DoFs are those in the closure of this set.
We refer to this method as~\ASMVankaStarPC{} (following~\cite{Rafiei_2024}) and illustrate it
with blue and red shaded areas in~\cref{fig:patches}\kern1pt\footnote{We acknowledge the use of tikz templates
from~\cite{farrell2021reynolds} in creating the mesh (\cref{fig:meshes}) and patch
(\cref{fig:patches}) figures.} for velocity and pressure DoF sets, respectively. 
The indices associated with these unknowns are then used to assemble the partitioning and local patch operators, which are solved exactly using LU factorization.  
Due to the discontinuous pressure field, the Scott-Vogelius discretization, unlike the Taylor-Hood, includes pressure DoFs on the boundary of the patch, which we denoted by red dashed squares to emphasize that these are ``internal'' degrees-of-freedom due to the DG discretization. 
For the FBF preconditioner, we apply patch-based relaxation only to the velocity block, defining patches based only on the vertex star (and not its closure, as needed for Vanka relaxation).  Following the naming convention in the Firedrake code~\cite{FiredrakeUserManual}, we call this the \ASMStarPC{} approach.

Firedrake~\cite{FiredrakeUserManual} provides two families of approaches for implementing patch relaxation methods, one based on using callbacks that allow matrix-free matrix-vector products for high-order operators but requires separately assembling the matrix for each patch~\cite{farrell2021pcpatch}, and one that assembles the global matrix and extracts the patch sub-matrices from this~\cite{Rafiei_2024}, see~\cref{rem:patch_setup}.  Here, we use the assembled matrix form, since it showed more efficient performance at large polynomial order (due to the significant overlap in patches, which leads to much duplication of effort in the callback approach).  Nonetheless, this explicit matrix assembly becomes a computational and storage disadvantage at high discretization orders, where one would prefer a matrix-free approach.
To alleviate this, \cite{brubeck2021scalable,brubeck2022multigrid} have introduced basis-change strategies that substantially improve complexity for relaxation on the patch systems associated with high-order discretization of spaces within the de Rham complex.
However, these approaches have not yet been extended to the triangular patches considered here, and applying these methods to coupled problems remains unexplored, highlighting an interesting area for future research.

Alternative relaxation schemes, such as Braess-Sarazin~\cite{Braess,he2019local}, Uzawa~\cite{gaspar2014simple,luo2017uzawa}, and distributive relaxation~\cite{wang2013multigrid,elman2008taxonomy} could also be considered in this setting. In addition to potentially requiring additional relaxation parameter tuning, these methods are not generally known to be robust with respect to discretization order. However, they could still be used on the coarser grids to achieve faster time to solution as demonstrated in~\cite{voronin2023monolithic}. To keep the focus on examining \hmg{} and \phmg{} schemes, the use of mixed relaxation schemes within a single preconditioner is not further explored in this work.

\begin{remark}[Comparison of Vertex-star Patch Setup Approaches]\label{rem:patch_setup}
As mentioned above, patch-based relaxation schemes can be defined in Firedrake using either \ASMStarPC{} and \ASMVankaStarPC{} or \PCPatch{}. 
The callback approach of \PCPatch{} demonstrates a notably quicker setup phase at lower discretization orders, when patches are small with fewer overlapping DoFs.
This approach, however, proves to be suboptimal at higher discretization levels where the overlap of DoFs between patches increases, resulting in substantial redundant computation.  
On the other hand, \ASMStarPC{} shows reasonable preformance at these higher orders, by effectively reducing the duplication of work associated with the overlap, thereby improving setup efficiency and performance.

We also note that \PCPatch{} is more difficult to extend, and currently lacks the functionality to construct the same patches as \ASMVankaStarPC{}.
This limitation precludes a balanced comparison between \PCPatch{} and \ASMVankaStarPC{} patches for the \afbf{} and \phmg{} solvers, respectively.
As a result, we opt out for a fairer comparison of \ASMStarPC{} and \ASMVankaStarPC{} patches for \afbf{} and \phmg{} solver, respectively. However, the effectiveness of \PCPatch{} for \afbf{} preconditioners at lower discretization orders underscores its value in specific computational contexts.
\Cref{table:patch_summary} provides a summary for the comparison of the \PCPatch{},
    \ASMVankaStarPC{}, and \ASMStarPC{} methods.
\end{remark}

\begin{table}[htbp]
\centering
\begin{tabular}{@{}lcp{8cm}@{}}
\toprule
\textbf{Patch Type}     & \textbf{Usage}  & \textbf{Characteristics} \\ \midrule
\PCPatch{}              & \phmg{}         & Faster setup at low discretization orders; suboptimal at high orders due to increased DoF overlap and redundant computation. Lacks the functionality to construct high-order monolithic Stokes patches. \\
\ASMVankaStarPC{}       & \phmg{}         & Effective at high discretizations using a global matrix to reduce overlap redundancy. Capable of constructing high-order monolithic Stokes patches. \\
\ASMStarPC{}            & \afbf{}         & Similar to \ASMVankaStarPC{}, but specifically used to assemble topological star patches such as vertex-star. \\ \bottomrule
\end{tabular}
\caption{Comparison of PCPatch, ASMVankaStarPC, and ASMStarPC.}
\label{table:patch_summary}
\end{table}

\section{Numerical Results}\label{sec:num_results}
\begin{figure}[!ht]
\begin{subfigure}{0.92\textwidth}
	\resizebox{\textwidth}{!}{
		\begin{tikzpicture}[%
		  every node/.style={draw=black, thick, anchor=west},
		  grow via three points={one child at (0.0,-0.8) and
		  two children at (0.0,-0.8) and (0.0,-1.6)},
		edge from parent path={(\tikzparentnode.210) |- (\tikzchildnode.west)}]
		\node {monolithic \protect\hmg{} V-cycle based on $K_\ell$}
		child {node {Relaxation on $K_\ell$, $\ell\in[0,\ell_{\text{max}})$}
		      child {node {Chebyshev ($\nu_h$ iterations)}
		         child {node {\protect\ASMVankaStarPC{} relaxation (1 iteration)}}
		     }
		     child[missing]{}
		}
		child[missing]{}
		child[missing]{}
		child {node {Grid transfer between $K_\ell$ and $K_{\ell+1}$}}
		child {node {Coarse grid solver}
		child {node {LU-factorization}}
		};
		\end{tikzpicture}
	}
	\caption{Schematic of the monolithic \protect\hmg{} V-cycles depicted in~\protect\cref{fig:hmg_cycle}.}
	\label{fig:hmg_precond}
\end{subfigure}   
\hfill 
\vspace{0.5cm}
\begin{subfigure}{0.92\textwidth}
	\resizebox{\textwidth}{!}{
		\begin{tikzpicture}[%
		  every node/.style={draw=black, thick, anchor=west},
		  grow via three points={one child at (0.0,-0.8) and
		  two children at (0.0,-0.8) and (0.0,-1.6)},
		edge from parent path={(\tikzparentnode.210) |- (\tikzchildnode.west)}]
		\node {monolithic \protect\phmg{} V-cycle ($G_\text{dc}$)}
		child {node {$p$\kern1pt-coarsening}
		   child {node {Relaxation on $K_0$}
		          child {node {Chebyshev ($\nu_p$ iterations)}
		             child {node {\ASMVankaStarPC{} relaxation (1 iteration)}}
		         }
		         child[missing]{}
		   }
		   child[missing]{}
		   child[missing]{}
		child {node {Grid transfer between $K_0$ and $K_1$}}
		}
	    child[missing]{}
	    child[missing]{}
	    child[missing]{}
	    child[missing]{}
		child {node {$h$MG based on $K_1$ ($G_1$)}
		};
		\end{tikzpicture}
	}
	\caption{Schematic of the monolithic \protect\phmg{} V-cycles depicted in~\protect\cref{fig:phmg_direct}.}
	\label{fig:phmg_precond}
\end{subfigure}
\hfill 
\vspace{0.5cm}
\begin{subfigure}{\textwidth}
	\resizebox{\textwidth}{!}{
		\begin{tikzpicture}[%
		  every node/.style={draw=black, thick, anchor=west},
		  grow via three points={one child at (0.0,-0.8) and
		  two children at (0.0,-0.8) and (0.0,-1.6)},
		  edge from parent path={(\tikzparentnode.210) |- (\tikzchildnode.west)}]
		  \node {Full-Block Factorization preconditioner}
		      child { node {Approximate Schur complement inverse}
		          child { node {Exact pressure mass matrix inverse}}
		      }
		      child [missing] {}
		      child { node {\protect\phmg{} V-cycle on momentum block}
					   child {node {Relaxation on $M_\ell$ for $\ell\in[0,\ell_{\text{max}})$}
					          child {node {Chebyshev ($\nu$ iterations)}
					             child {node {ASMStarPC relaxation (1 iteration)}}
					         }
					         child[missing]{}
					   }
		};
		\end{tikzpicture}
	}
	\caption{Schematic of the \protect\afbf{}{} block preconditioner given by~\protect\Cref{eq:fbf}. }
	\label{fig:afbf_solver}
\end{subfigure}
\vskip -.05in
\caption{Overview of the high-order preconditioners used for solving~\protect\cref{eq:saddle}. }     
\label{fig:phmg_combined}
\end{figure}
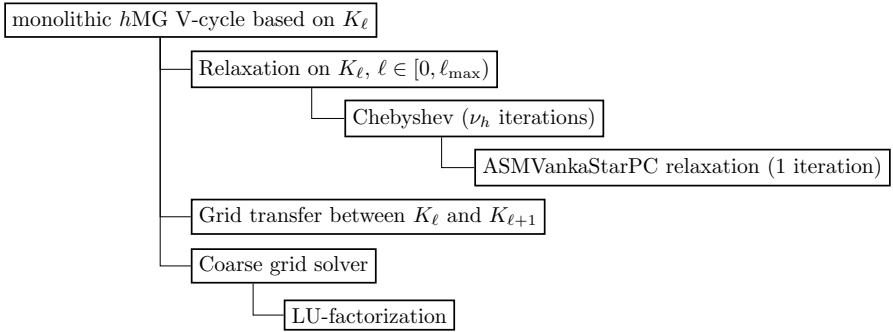
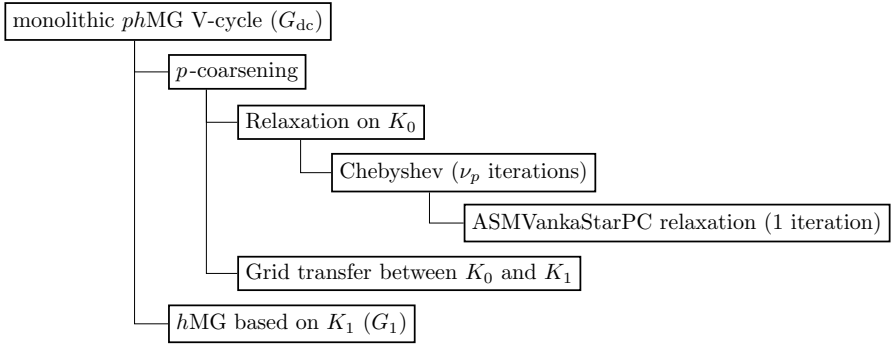
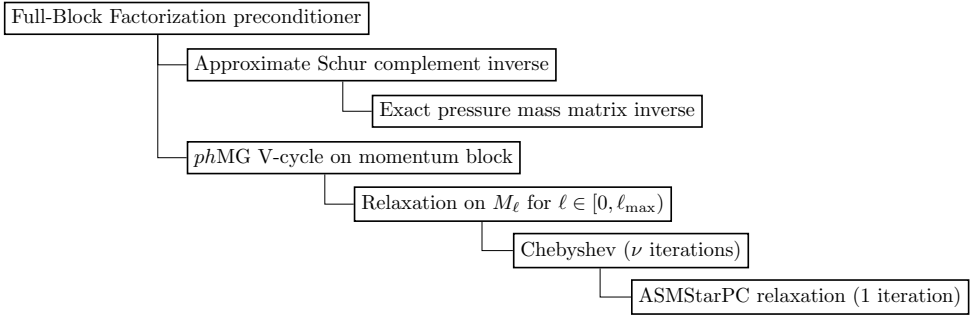

\begin{figure}[htbp]
\resizebox{1\textwidth}{!}{
\begin{tikzpicture}[>=Stealth,font=\Large]

  \newcommand{\hunit}{-2.0}
  \newcommand{\vdist}{-1.0}

  \newcommand{\vstart}{0}

  \tikzset{Pnode/.style={rectangle, draw, thin, minimum size=2.1em, inner sep=0}}
  
   \draw [decorate,decoration={brace,amplitude=7pt,mirror,aspect=0.5},yshift=-3pt,thick] 
    (-5,5.7) -- (-5,10.5) node [black,midway,yshift=-0.0cm,xshift=1.5cm] {3 $p$-grids};

   \draw [decorate,decoration={brace,amplitude=7pt,mirror,aspect=0.5},yshift=-3pt,thick] 
    (-5,2.6) -- (-5,5.5) node [black,midway,yshift=-0.0cm,xshift=1.5cm] {2 $p$-grids};

  \foreach \x [evaluate=\x as \hpos using \x*\hunit+0.5*\hunit] in {10,9,...,3}{
    \node[Pnode] (P\x) at (\hpos-1,\vstart-\x*\vdist) {\( \pmb{\mathbb{P}}_{\x} \)};

    \pgfmathtruncatemacro{\target}{int(\x==5||\x==4||\x==3 ? 2 : (\x<8 ? 4 : 5))}

    \ifnum\target=2
      \node[Pnode] (P\target) at (3*\hunit,\vstart-\x*\vdist) {\( \pmb{\mathbb{P}}_{\target} \)};
      \draw[dashed,->] (P\x) -- (P\target);
    \else
      \node[Pnode] (P\target\x) at (\target*\hunit+\hunit,
                                    \vstart-\x*\vdist) {\( \pmb{\mathbb{P}}_{\target} \)};
      \draw[dashed,->] (P\x) -- (P\target\x);
      \node[Pnode] (P2\x) at (3*\hunit,\vstart-\x*\vdist) {\( \pmb{\mathbb{P}}_2 \)};
      \draw[dashed,->] (P\target\x) -- (P2\x);
    \fi
  }
\end{tikzpicture}
}
\caption{$p$MG coarsening schedule for velocity field. For \protect\vecP{k} with $k\in[2,5]$, coarsening maps directly to \protect\vecP{2}. For  
$k\in[6,7]$, we first coarsened to \protect\vecP{4} and then to \protect\vecP{2}. For 
$k\in[8,10]$, we initially coarsen to \protect\vecP{5}, followed by \protect\vecP{2}.}
\label{fig:pmg_schedule}
\end{figure}
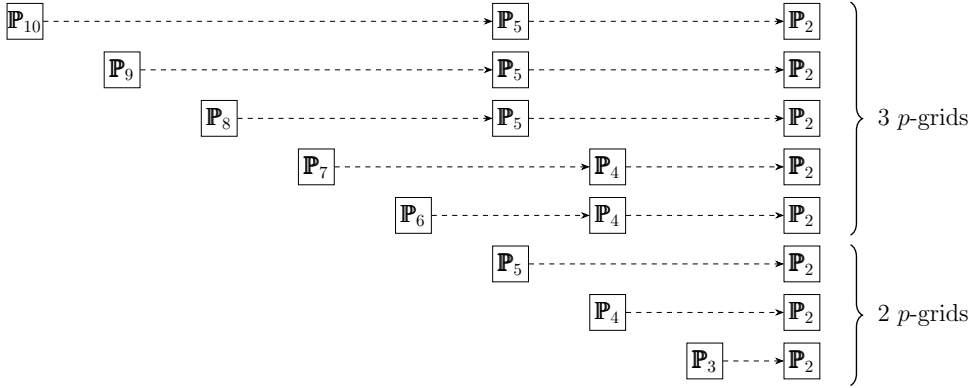

Numerical experiments are conducted in the Firedrake
framework~\cite{FiredrakeUserManual}, which uses
PETSc~\cite{balay2019petsc} and includes patch-based relaxation schemes through PCPATCH~\cite{farrell2021pcpatch} and \ASMVankaStarPC{}~\cite{Rafiei_2024}.  ALFI\kern1pt\footnote{ALFI codebase \href{https://github.com/florianwechsung/alfi/}{https://github.com/florianwechsung/alfi/}.} is also used to facilitate barycentric mesh refinement needed for \SVPk{}.  \Cref{fig:phmg_combined} provides solver diagrams for the preconditioners introduced in the prior section.
Coarse-grid operators are computed via rediscretization, which is equivalent to
Galerkin coarsening with our grid transfer choices.  For \phmg{}, the spectral
coarsening of the velocity field always concludes with the~\vecP{2} basis.  For
monolithic \phmg{}, the pressure basis remains one order below the velocity
basis, giving rise to either \THPk{} or \SVPk{} discretizations. Our results
consider aggressive and gradual $p$\kern1pt-coarsenings.  \VphmgDirect{}
performs $p$\kern1pt-coarsening in one step directly to \THP{2}.
\Cref{fig:pmg_schedule} illustrates the coarsening schedule for \VphmgGradual{}.
Damping parameters within relaxation are effectively defined through Chebyshev
iteration, which uses an estimate for the largest eigenvalue of $M_\ell^{-1}
K_\ell$.  On the coarsest grid, MUMPS is used to perform a direct solve.
Code and implementation details can be found at \href{\githuburl}{\githuburltext}.

In the following subsection, we present performance comparisons of \phmg{} solvers versus \hmg{}
using restarted FGMRES with a maximum Krylov subspace size of $30$, where convergence is declared when 
the the $\ell_2$-norm of the residual is reduced by a relative factor of $10^{10}$. 
We focus on two Stokes flow problems on triangular (tetrahedral) meshes:
\begin{itemize}
\item regularized lid-driven cavity over structured meshes; and
\item backward-facing step over unstructured meshes.
\end{itemize}
These problems are detailed in Examples 3.1.2 and 3.1.3 of~\cite{elman2014finite}. 
Our reporting highlights multigrid parameter sets yielding the fastest average convergence for the entire set of discretizations and refinement levels rather than the best solver for individual discretization and mesh refinement. Consequently, not all combinations are discussed.

\subsection{Taylor-Hood}\label{sec:th}

This section presents results for \THPk{} discretizations. 

\subsubsection{2D Case}\label{sec:th_2D}

\begin{figure}[!ht]
\begin{subfigure}{1\textwidth}
\includegraphics[width=\textwidth]{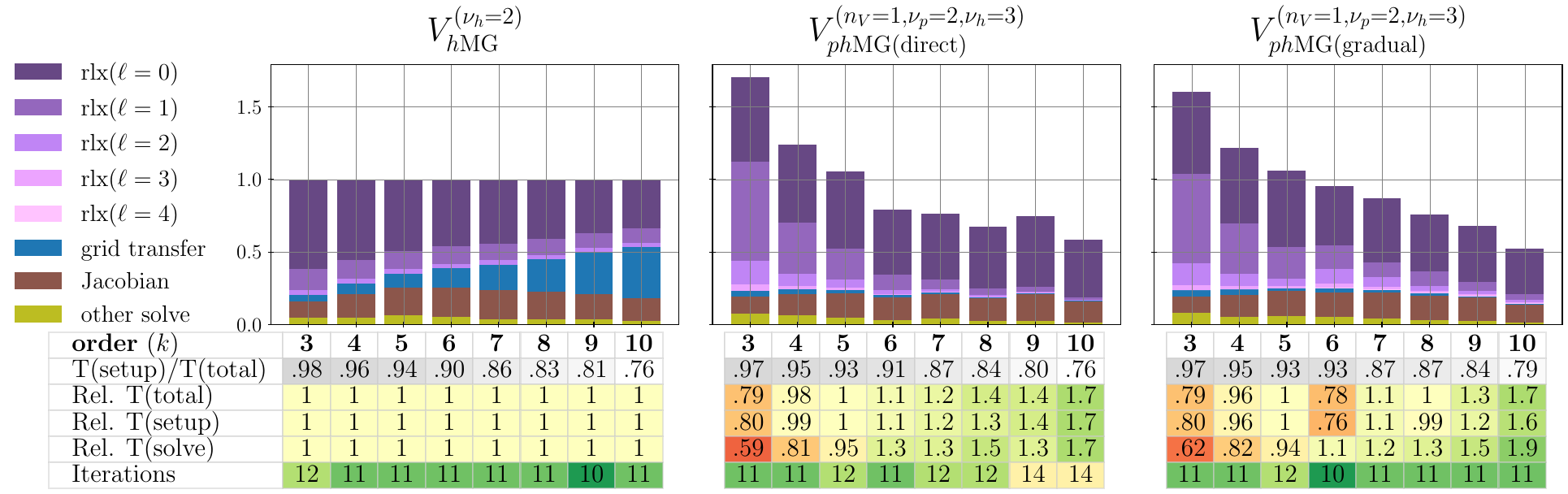}
\caption{Lid-driven cavity problem on structured mesh.}
\label{fig:th_ldc2d_ref3}
\end{subfigure}
\begin{subfigure}{1\textwidth}
\includegraphics[width=\textwidth]{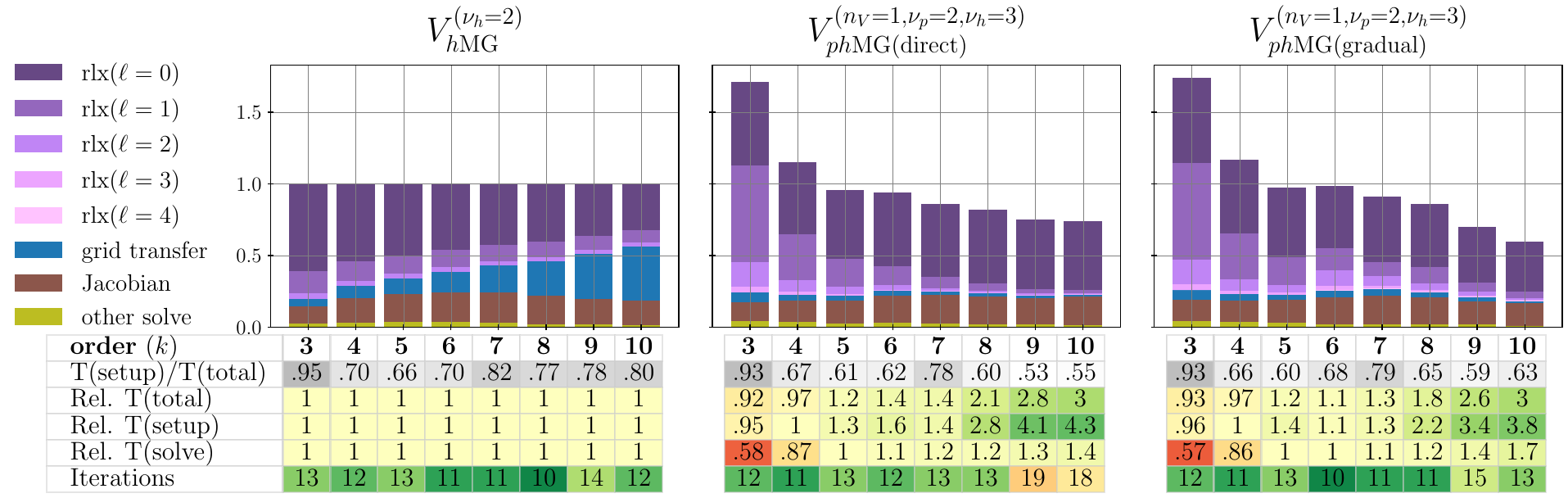}
\caption{Backward facing step on unstructured mesh.}
\label{fig:th_bfs2d_ref3}
\end{subfigure}
\caption{2D Taylor-Hood (\THPk{}): Deep hierarchies with 3 mesh refinements. The full solver name is abbreviated in text with \hmg{}, \VphmgDirect{}, and \VphmgGradual{}, respectively.}
\label{fig:th_2D_ref_3}
\end{figure}

This section presents results obtained on an Intel Xeon CPU E5-2650 v3 @ 2.30GHz with 256GB RAM, utilizing 8 MPI tasks to balance computation/communication trade-offs. 
Memory limitations required the use of smaller coarser grids for problems involving higher-order computations due to the increased stencil size. The number of unknowns and the average number of non-zeros per DoF in the finest-grid matrix operator for each problem shown in~\cref{fig:th_2D_ref_3} are listed in~\cref{table:nnz_th_2d}.
\begin{table}%
\centering
\begin{tabular}{c|c|rrrrrrrr}
\toprule
& & \multicolumn{8}{c}{\textbf{order} ($k$)} \\
\textbf{problem type} & \textbf{metric} & 3 & 4 & 5 & 6 & 7 & 8 & 9 & 10 \\
\midrule
\multirow{2}{*}{ldc2d} & DoFs (mil.) & 3.7 & 4.0 & 4.1 & 4.2 & 4.2 & 4.5 & 4.7 & 4.6 \\
                       & nnz/DoFs   & 45 & 63 & 85 & 109 & 137 & 167 & 201 & 237 \\
\midrule
\multirow{2}{*}{bfs2d} & DoFs (mil.) & 2.5 & 2.4 & 2.0 & 3.0 & 2.2 & 1.6 & 1.0 & 1.3 \\
                       & nnz/DoFs   & 45 & 63 & 85 & 109 & 137 & 167 & 200 & 237 \\
\bottomrule
\end{tabular}
\caption{Total number of degrees of freedom (DoFs, in millions) and the average number of non-zero entries per DoF in the Stokes operator discretized using the Taylor-Hood element pair. The data is presented for different approximation orders ($k$) across two problems: the 2D lid-driven cavity (ldc2D) and the backward-facing step (bfs2D).}
\label{table:nnz_th_2d}
\end{table}

The bar plots in~\cref{fig:th_2D_ref_3} illustrate the distribution of time spent during the solution phase across the main kernels, normalized against the total solve time of the \hmg{} solver for the same polynomial order for comparison. 
Below each bar plot, a data table presents key metrics, where total time is defined as the combined duration of the solver setup (including low-order rediscretization) and the solve phase:
\begin{center}%
    \begin{tabular}{p{3.0cm}|p{2.2cm}|p{6.5cm}}
\toprule
\textbf{Metric} & \overfullrule=0pt\textbf{Abbreviation} &  \textbf{Description} \\
\midrule
\overfullrule=0pt\setupOverTotalFull{} & \setupOverTotal{} & Setup phase time as a fraction of total time \\
\relTotalFull{} & \relTotal{} & Total time relative to the reference solver \\
\relSetupFull{} & \relSetup{} & Setup time relative to the reference solver \\
\relSolveFull{} & \relSolve{} & Solve time compared to the reference solver \\
\textbf{Iteration} & \textbf{iterations} & Number of FGMRES iterations. \\
\bottomrule
\end{tabular}
\label{table:metrics}
\end{center}
The relative metrics indicate the speed-up compared to the reference solvers. For example, \relTotal{} metric for \phmg{} solvers is computed as follows:
$$\textbf{R(total)}= \frac{\textbf{Total time (\hmg{})}}{\textbf{Total time (\phmg{})}},$$
where \hmg{} is the reference solver. 

The \setupOverTotal{} metric is visualized using a grey-to-white spectrum. Darker shades signify a larger fraction of setup time.
Subsequent relative metrics are highlighted in shades of {\color{tab:green}green} or {\color{tab:red}red}, where green indicates faster and red denotes slower relative to the fastest solver.  We also highlight iteration counts with lower values in shades of {\color{tab:green}green} and higher iteration counts in shades of {\color{tab:red}red}.

The convergence and timings results in~\cref{fig:th_2D_ref_3} document the data for \hmg{} and \phmg{} solvers with 4 levels of $h$-grids. \phmg{} includes 1-2 additional $p$\kern1pt-refinement grids, which vary based on the coarsening type (either direct or gradual) and the \VphmgGradual{} coarsening schedule, as illustrated in~\cref{fig:pmg_schedule}. \Cref{fig:th_2D_ref_3} presents results for the 2D Taylor-Hood (\THPk{}) problems, including the lid-driven cavity on a structured mesh and the backward-facing step on an unstructured mesh.

The setup phase ratios in~\cref{fig:th_2D_ref_3} indicate 
a significant portion of the total time, ranging from 53\% to 98\% is consumed during the setup phase of the preconditioner.
Within each problem type and discretization order, all solvers show approximately similar setup phase times. Across problem types, we observe that the backward-facing step has a lower setup-phase ratio than the same solvers in lid-driven-cavity problems. 

The \relTotal{} performance of \phmg{} solvers compared to \hmg{} varies with discretization order; they are slower or the same speed for $k\in[3,5]$, but show improvement at higher orders, where \phmg{} achieves speedups of up to 1.7x for the lid-driven cavity and up to 3x for the backward-facing step. 
The observed speed-ups for the backward-facing step are primarily due to a more efficient setup phase, which benefits from sparser coarse grids and less costly coarse-grid relaxation. The added overhead cost of $p$\kern1pt-coarsening at lower discretization orders makes the \phmg{} less computationally efficient than \hmg{} in both the setup and solve phases.
This inefficiency arises for several reasons. Firstly, the \phmg{} cycles involve more levels. Moreover, when performing $p$-coarsening, systems are assembled on the fine mesh, and even for low-order systems, those operators are costly. For instance, as can be observed from the bar plots for \VphmgDirect{} solvers, a significant amount of time is spent on the \THP{2} relaxation ($\ell=1$). However, as the discretization order $k$ increases, the relative cost of relaxation on $\ell=1$ becomes less significant compared to $\ell=0$. Similar logic applies to \VphmgGradual{} solvers.

 The speed-up of the \phmg{} setup is primarily attributed to sparser coarse grids and cheaper relaxation setup. Secondary factors include less expensive grid-transfer costs. The coarse-grid solve assembly amounts to an insignificant part of the overall run-time for these 2D problems.   
During the solve phase, the bar plots show in~\cref{fig:th_2D_ref_3} that for lower values of $k$, the solve phase in \hmg{} solvers is dominated by relaxation costs. However, as $k$ increases, the grid-transfer costs increase from less than 10\% to approximately one-third of the solve phase.  
\phmg{} avoids this growth in interpolation cost by using lower-order bases despite having more levels in the multigrid due to $p$\kern1pt-coarsening. At low $k$ values, \phmg{} is rarely faster than \hmg{}. 
However, as the discretization order increases, \phmg{} shows improved relative solve phase timings, benefiting from less expensive grid transfer and less costly relaxation on lower-order coarse grids.

Shifting attention to comparing the \VphmgDirect{} and \VphmgGradual{} solvers, the hierarchies of \VphmgDirect{} and \VphmgGradual{} coincide for $k\in[3, 5]$, due to identical $p$\kern1pt-coarsening schedules.
However, for $k\geq 6$, \VphmgGradual{} solvers introduce an additional $p$\kern1pt-coarsening grid. \sloppy This gradual $p$\kern1pt-coarsening allows \VphmgGradual{} solvers to maintain convergence rates similar to \hmg{}, whereas \VphmgDirect{} solvers experience an increase in iteration counts with rising $k$, particularly evident in unstructured base meshes~\cref{fig:th_bfs2d_ref3,fig:th_bfs2d_ref3}.

As described at the beginning of~\cref{sec:num_results}, we have explored a variety of relaxation schedules as well as different numbers of \hmg{} V-cycles in \phmg{} solvers.
The reported results showcase the fastest combination of $(n_V, \nu_p,\nu_h)$
parameters, yet other configurations, notably \Vphmg{2}{2}{2}{direct} and
\Vphmg{2}{2}{2}{gradual}, which execute two \hmg{} V-cycles per iteration, also
display competitive performance. A common feature among these \phmg{} solvers is
the need for an increased accuracy of \hmg{} corrections, the lower part of the \phmg{} cycle, to maintain low iteration counts.
An alternative method to achieve this efficacy involves one-sided V-cycles
utilizing high-order Chebyshev polynomials, which is a robust strategy for
managing poor-quality coarse spaces due to anisotropy or aggressive
coarsening~\cite{phillips2022optimal}.

\begin{remark}[Bias in the Convergence and Timings Results] \label{rmrk:bias}
A notable pattern is observed with all \phmg{} solvers: they tend to outpace \hmg{} solvers more significantly at higher discretization orders. 
This trend is partially due to the selection of a solver configuration optimized for the fastest overall time to solution across a range of discretization orders ($k\in[3,10]$). Since higher-order problems generally require more time for setup and solution phase, this parameter optimization approach favors them. 
More suitable parameter combinations (specifically, $n_V, \nu_p, \nu_h$) could be identified for lower-order discretizations to improve efficiency. 
However, the exploration of these optimal configurations for each discretization level, considering the constraints of time, computing resources, and the diversity of problems, is beyond the scope of this study.
\end{remark}

\begin{remark}[\phmg{} Relaxation at Low Discretization Order]
In examining the bar plots for \phmg{} solvers at low discretization orders,
    particularly $k=3$ and $k=4$, it is observed that the time spent on
    relaxation at levels $\ell=0$ and $\ell=1$ is quite similar, with relaxation
    at $\ell=1$ sometimes exceeding that at $\ell=0$. This pattern is due to the
    \phmg{} configuration conducting more relaxation sweeps within the \hmg{}
    hierarchy compared to $p$\kern1pt-coarsening levels. Despite this, relaxation at the coarser levels contributes significantly to the total solve phase time. 
In contexts of low-order discretization, block-wise relaxation methods like the least-squares-commutator~\cite{wang2013multigrid} may be more computationally efficient than patch relaxation, provided that convergence is not compromised.
\end{remark}

\subsubsection{3D Problems}

This section presents the numerical results of 3D Taylor-Hood problems obtained using the NCSA Delta supercomputer\kern1pt\footnote{NCSA Delta \hyperlink{https://www.ncsa.illinois.edu/research/project-highlights/delta/}{https://www.ncsa.illinois.edu/research/project-highlights/delta/}.}. Each computing node consists of two AMD EPYC 7763 ``Milan'' CPUs, each with 128 cores, operating at a frequency of 2.45 GHz, and equipped with 256 GB of shared RAM.

\begin{figure}%
\centering
\begin{subfigure}{1\textwidth}
\includegraphics[width=\textwidth]{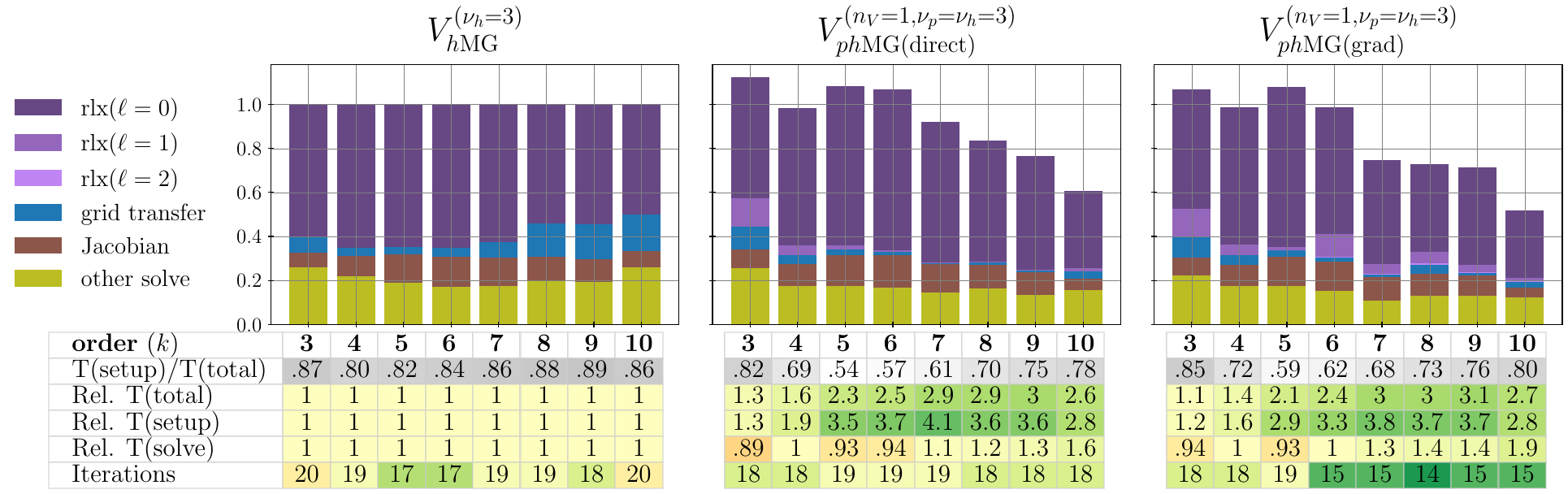}
\caption{Results for hierarchies with single mesh refinement.}
\label{fig:th_ldc3d_ref1}
\end{subfigure}
\begin{subfigure}{1\textwidth}
\includegraphics[width=\textwidth]{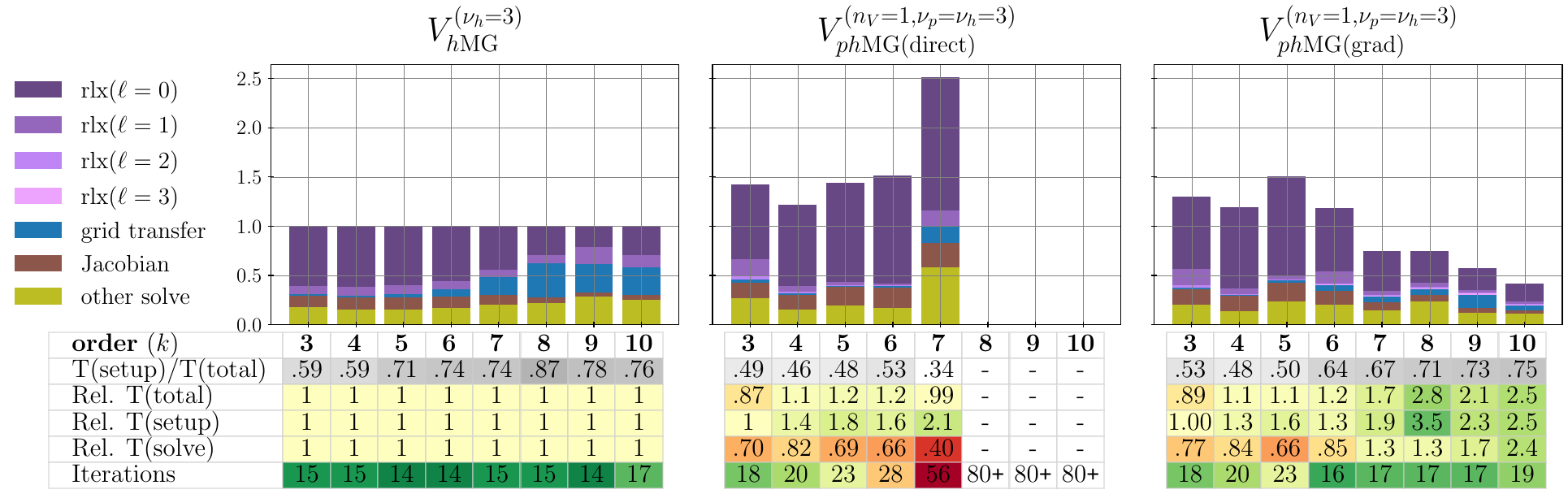}
\caption{Results for hierarchies with two mesh refinements.}
\label{fig:th_ldc3d_ref2}
\end{subfigure}
\caption{3D Taylor-Hood (\THPk{}) Lid-driven cavity on structured mesh.}
\label{fig:th_ldc3d_ldc}
\end{figure}

We focus on assessing the efficacy of the \hmg{} and \phmg{} solvers when addressing the lid-driven-cavity problem on cubic structured meshes. Our analysis includes results from employing one and two $h$-refinements, utilizing 72 and 192 MPI tasks, respectively.
Similar to the 2D scenario, we formulated higher-order problems on smaller base meshes to optimize memory usage.
The number of unknowns and the average number of non-zero entries per degree of freedom (DoF) in the matrix operator for each problem shown in~\cref{fig:th_ldc3d_ref2} are listed in~\cref{table:nnz_th_3d}. As can be seen in~\cref{table:nnz_th_3d} the number of non-zeros values in the matrix grows much faster for 3D than 2D problems. 
\overfullrule=0pt\begin{table}[htbp]
\centering
\overfullrule=0pt\begin{tabular}{p{2.4cm}|p{1.8cm}|rrrrrrrr}
\toprule
\multicolumn{2}{c|}{} & \multicolumn{8}{c}{\textbf{order} ($k$)} \\
\textbf{problem type} & \textbf{metric} & 3 & 4 & 5 & 6 & 7 & 8 & 9 & 10 \\
\midrule
\multirow{2}{*}{ldc3d} & DoFs (mil.) & 15 & 18 & 20 & 25 & 17 & 15 & 11 & 15 \\
                       & nnz/DoFs   & 171 & 270 & 400 & 565 & 766 & 1010 & 1298 & 1641 \\
\bottomrule
\end{tabular}
\caption{Total number of degrees of freedom (DoFs, in millions) and the average number of non-zero entries per DoF in the Stokes operator discretized using the Taylor-Hood element pair. The data is presented for different approximation orders ($k$) for 3D lid-driven cavity (ldc3D) problem with two mesh refinements (\cref{fig:th_ldc3d_ref2}).}
\label{table:nnz_th_3d}
\end{table}

It is important to address the significant memory requirements of 3D monolithic multigrid solvers using patch relaxation. 
These constraints necessitate a specific configuration of computational resources. The choice of the number of MPI tasks, which is not a straightforward multiple of the cores available per node, was heuristically identified to accommodate the memory requirements of the highest-order discretizations. 
This approach balances the constraints of memory capacity with the goal of
ensuring effective workload distribution across tasks. In dealing with 3D
problems, our aim is to have a minimum of $2\texttt{+}$ elements per MPI task on
the coarsest grid for high-order discretizations, scaling to hundreds of
elements per MPI task at lower-order discretizations.

 Similar to the 2D Taylor-Hood findings, we observe in~\cref{fig:th_ldc3d_ref1,fig:th_ldc3d_ref2} that the setup phase accounts for a large fraction of total time, ranging from 59\% to 89\% for \hmg{} solvers and 34\% to 82\% for \phmg{} solvers. The fraction of time spent in the setup phase appears to decrease slightly when going from solvers with 1 to 2 mesh refinements. However, the overall change does not appear to be significant. 

All solvers have approximately stable iteration counts when examining the iteration counts within shallow hierarchies, as shown in~\cref{fig:th_ldc3d_ref1}. This stability is largely attributed to the effectiveness of exact coarse-grid correction.
However, challenges arise in solvers with deeper hierarchies, specifically those with two mesh refinements (\cref{fig:th_ldc3d_ref2}). 
In these cases, \VphmgDirect{} struggles to converge in under 80 iterations starting from $k=7$. Although the iteration count can be somewhat reduced by increasing the number of \hmg{} V-cycles (not depicted), the issue of growing iteration counts remains. 
In contrast, the \VphmgGradual{} solver demonstrates consistent iteration counts by integrating intermediate $p$\kern1pt-coarsening grids.

Unlike with 2D problems, \phmg{} solvers in 3D begin to outperform \hmg{} solvers in terms of total relative time (\relTotal{}) starting with lower $k$ values.
This advantage is primarily attributed to faster-growing linear-operator stencils for high-order operators in 3D, making the setup and solution phases significantly more computationally efficient in \phmg{} solvers compared to \hmg{} solvers. This leads to improved relative setup timings (\relSetup{}) and a more pronounced difference in computational cost during the setup and solve phases in 3D problems.

During the solve phase, the \relSolve{} timings for \phmg{} solvers in 3D also reveal an improvement compared to their performance in 2D scenarios. 
In 2D cases, the most notable speed-ups in the solve phase were recorded at $k=10$, achieving 1.7x to 1.9x speedup for \VphmgGradual{} hierarchies. In the context of 3D problems, \phmg{} solvers reach speed-ups of 1.9x and 2.4x for the same discretization order ($k=10$). 
This observed improvement is partly due to the increased computational demands of coarse-level relaxation (rlx($\ell=1$)) and coarse-grid solving (grouped into ``other solve'' label) in 3D \hmg{} solvers, as indicated in the bar plots of~\cref{fig:th_ldc3d_ldc}.

\subsection{Scott-Vogelius}\label{sec:sv}

This section examines the convergence of multilevel preconditioners for Stokes problems discretized with Scott-Vogelius elements on barycentric meshes. 
Given the absence of a nested hierarchical multigrid (\hmg) preconditioner tailored for these meshes, we replace the monolithic \hmg{} reference solver with a full-block-factorization solver, \afbf{}, as described in~\cref{fig:afbf_solver}.
The same computing environments and the number of MPI tasks have been used as in Taylor-Hood 2D and 3D results, respectively.

\subsubsection{2D Problems}

\begin{figure}[!ht]
\begin{subfigure}{1\textwidth}
\includegraphics[width=\textwidth]{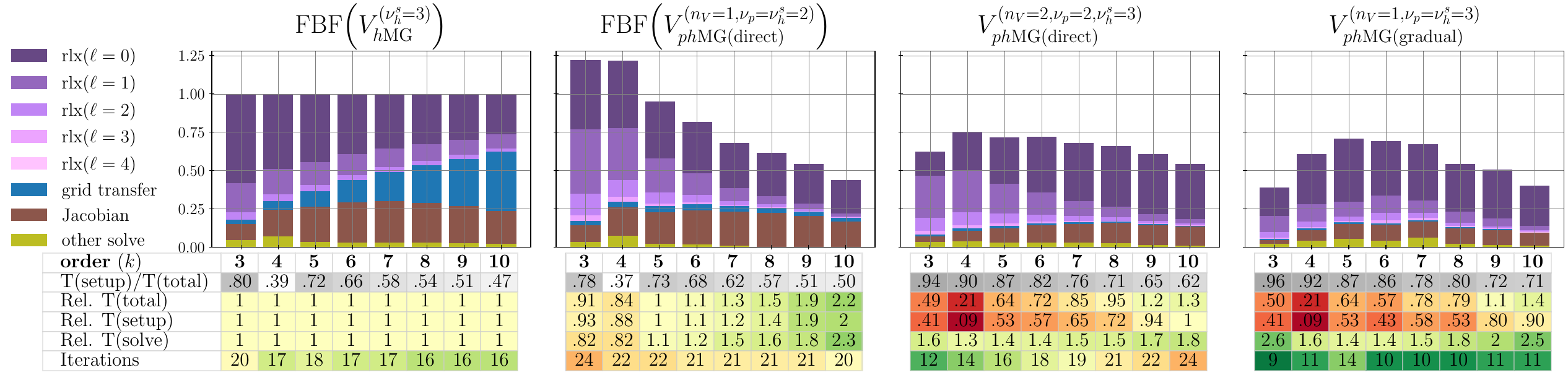}
\caption{Lid-driven cavity problem on structured mesh.}
\label{fig:sv_ldc2d_ref3}
\end{subfigure}
\begin{subfigure}{1\textwidth}
\includegraphics[width=\textwidth]{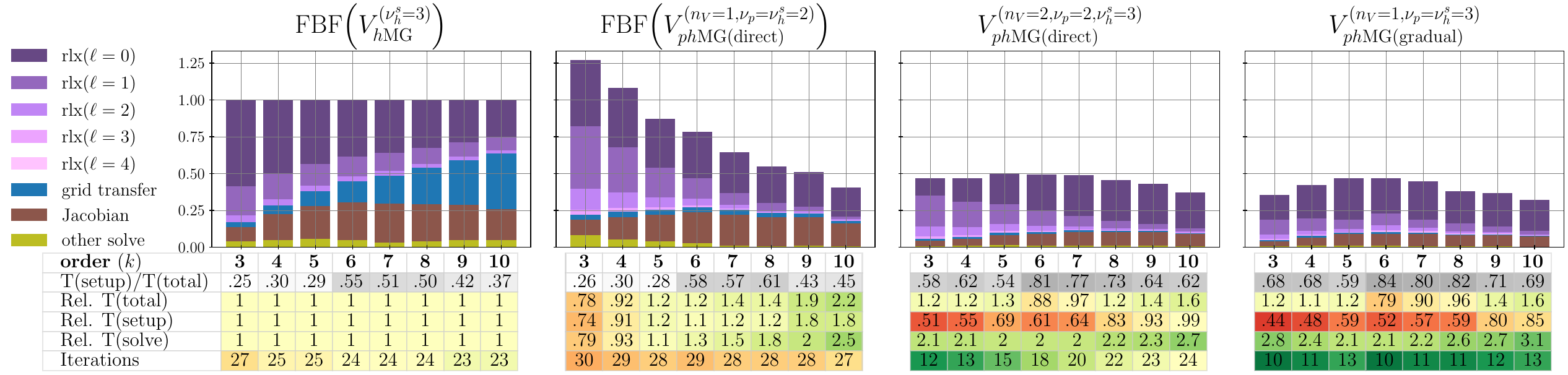}
\caption{Backward facing step on unstructured mesh.}
\label{fig:sv_bfs2d_ref3}
\end{subfigure}
\caption{2D Scott-Vogelius (\SVPk{}): Deep hierarchies with 3 mesh refinement. The full solver name is abbreviated in text with \fhmg{}, \fphmg{}, \VphmgDirect{} and \VphmgGradual{}, respectively.}
\label{fig:sv_2D}
\end{figure}

In~\cref{fig:sv_2D}, we present the convergence timing results for four distinct
solvers. Among these, two belong to the category of \afbf{} solvers. One of
these solvers employs \hmg{} for the velocity field, referred to as \fhmg{},
while the other utilizes \phmg{} with a direct coarsening approach on the
velocity field, denoted as \fphmg{}. The other two solvers are \VphmgDirect{}
and \VphmgGradual{}. 
The number of unknowns and the average number of non-zero entries per degree of freedom (DoF) in the matrix operator for each problem shown in~\cref{fig:sv_2D} are listed in~\cref{table:nnz_sv_2d}.
\begin{table}[htbp]
\centering
\begin{tabular}{p{2.4cm}|p{1.8cm}|rrrrrrrr}
\toprule
& & \multicolumn{8}{c}{\textbf{order} ($k$)} \\
\textbf{problem type} & \textbf{metric} & 3 & 4 & 5 & 6 & 7 & 8 & 9 & 10 \\
\midrule
\multirow{2}{*}{ldc2d} & DoFs (mil.) & 3.9 & 4.0 & 3.9 & 3.7 & 3.6 & 3.8 & 3.1 & 3.8 \\
                       & nnz/DoFs   & 38 & 55 & 75 & 98 & 124 & 153 & 185 & 220 \\
\midrule
\multirow{2}{*}{bfs2d} & DoFs (mil.) & 2.7 & 4.7 & 3.9 & 3.0 & 4.1 & 5.3 & 3.4 & 4.2 \\
                       & nnz/DoFs   & 38 & 55 & 75 & 98 & 124 & 153 & 185 & 220 \\
\bottomrule
\end{tabular}
\caption{Total number of degrees of freedom (DoFs, in millions) and the average number of non-zero entries per DoF in the Stokes operator discretized using the Scott-Vogelius element pair. The data is presented for different approximation orders ($k$) across two problems: the 2D lid-driven cavity (ldc2D) and the backward-facing step (bfs2D).}
\label{table:nnz_sv_2d}
\end{table}

A comparative analysis reveals a contrast in the setup time ratios (\setupOverTotal{}) between \afbf{} and \phmg{} solvers. The \afbf{} solvers exhibit consistently lower setup times, ranging from 25\% to 80\%, whereas \phmg{} solvers mirror the setup time seen in 2D Taylor-Hood analyses, with setup phase consuming a larger share of total runtime, ranging from 54\% to 96\%. Notably, for the backward-facing step problems on unstructured meshes, \setupOverTotal{} values for \phmg{} are not as large as they are for lid-driven cavity scenarios, suggesting a mesh-dependent behavior in the patch assembly process between \ASMStarPC{} in \afbf{} and \ASMVankaStarPC{} in \phmg{}.
The \relTotal{} metric shows that \fphmg{} significantly outperforms \fhmg{}, starting with $k=5$, achieving up to a 2.2x reduction in total solution time. 
The improved \relTotal{} of the \fphmg{} solver relative to \fhmg{} is due to consistently faster setup and solve phase times. Since most of the run-time is spent in the solution phase, the overall speed-up of \fphmg{} is closer to the \relSolve{} ratio. 

The \phmg{} solvers outperform \fhmg{} in terms of \relTotal{} only at the highest discretization orders, $k=9$ and $k=10$, for lid-driven cavity problems. 
For the backward-facing step problems, \phmg{} performs almost as well or better than \fhmg{} for both low and high discretization orders, but not around $k=6$ and $k=7$. 
\fphmg{} achieves significantly faster time to solution that \phmg{} for most problem sizes.
The slower relative performance of \phmg{} solvers with respect to \fphmg{} is attributed to much slower setup times due to the large cost associated with the setup of patches containing multiple unknown types.  

We note that, as the discretization order increases, both \fhmg{} and \fphmg{} show a reduction in the number of iterations required. 
While \fphmg{} takes slightly more iterations to converge than \fhmg{}, it achieves faster solve phase timings, \relSolve{}, than \hmg{} starting from $k=5$.
The improved performance of \fphmg{} at higher discretization order is attributed to its lower per-iteration cost due to the use of low-order coarse grids, which translates into cheaper coarse-grid relaxation, which is readily observed in the bar plots (e.g., rlx($\ell=1$)). 
In addition,  it is observed that \fphmg{} grid-transfer timings do not grow with $k$ as they do for \fhmg{}, a behavior that was also observed in \hmg{} Taylor-Hood problems.

Conversely, \VphmgDirect{} displays a pronounced increase in iteration counts as $k$ increases, more so than what is observed for Taylor-Hood problems. 
More rapid iteration growth necessitates additional \hmg{} cycles or relaxation sweeps relative to solvers presented in the Taylor-Hood section. 
Despite the increased iteration count, \VphmgDirect{} achieves a \relSolve{} performance that is significantly faster, up to 1.8 to 2.7 times, than \fhmg{}. It also shows a marginal improvement in solve time over \fphmg{} for the backward-facing step and a slight underperformance for the lid-driven-cavity problem.
\VphmgGradual{} mitigates the increase in iteration counts seen in \VphmgDirect{}, achieving significantly faster relative solve phase timings and outperforming all other solvers in \relSolve{} for most discretization orders.

\subsubsection{3D Problems}

\begin{figure}[!ht]
\centering
\begin{subfigure}{1\textwidth}
\includegraphics[width=\textwidth]{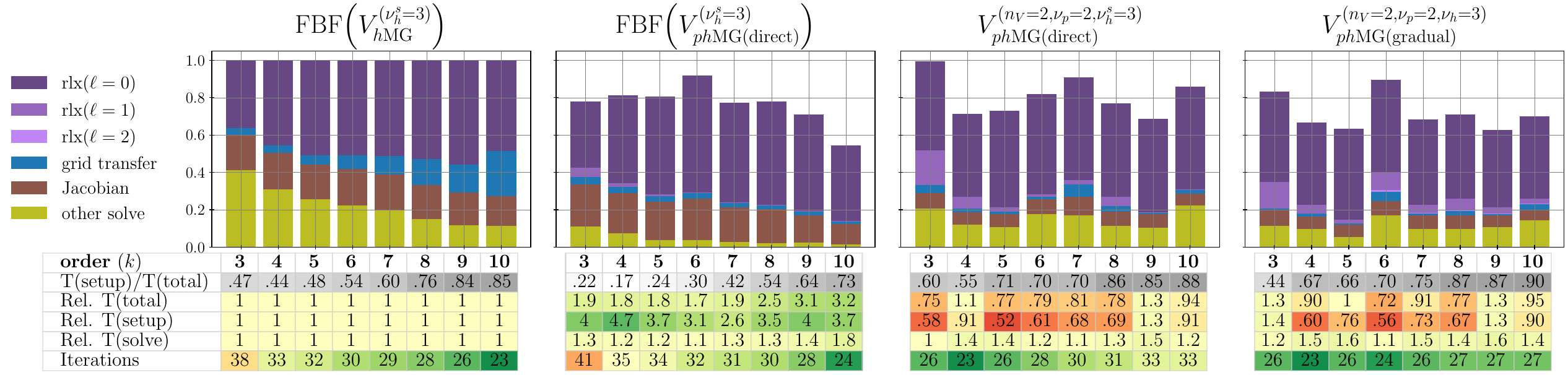}
\caption{Results for hierarchies with single mesh refinement.}
\label{fig:sv_ref1_ldc3d}
\end{subfigure}
\begin{subfigure}{1\textwidth}
\includegraphics[width=\textwidth]{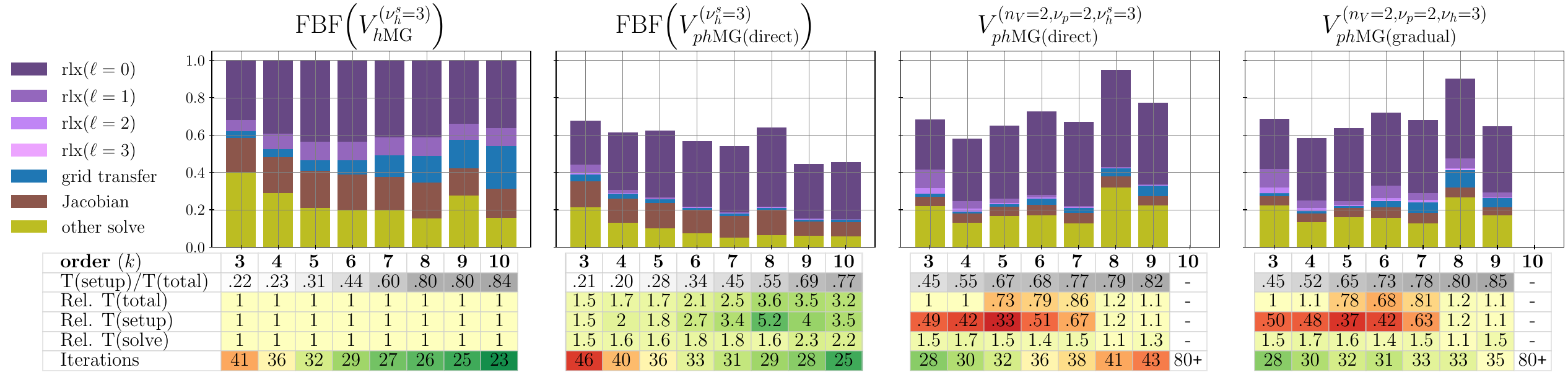}
\caption{Results for hierarchies with two mesh refinements.}
\label{fig:sv_ref2_ldc3d}
\end{subfigure}
\caption{3D Scott Vogelius(\SVPk{}) Lid-driven cavity on structured mesh.}
\label{fig:sv_ldc_3d}
\end{figure}

 \cref{fig:sv_ldc_3d} shows a noticeable deviation in the \setupOverTotal{} ratio trends for \afbf{} solvers compared to the 2D scenarios
as $k$ increases, with both \afbf{} preconditioners displaying an upwards trend in the fraction of overall time spent on the setup phase, going from about 20\% to 80\% or more. This increase is attributed to the growth of velocity patch sizes and the increase in the block sizes in the discontinuous pressure mass matrix as the discretization order increases.
Similarly, \phmg{} solvers show a steep increase in the \setupOverTotal{} ratio, more pronounced than in 2D Scott-Vogelius and 3D Taylor-Hood scenarios. This increase is primarily due to the expensive Scott-Vogelius vertex-star patches, encompassing a growing number of velocity and discontinuous pressure DoFs.  Due to the lack of convergence of \phmg{} solvers for $k=10$ problems at two mesh refinements, we leave those columns blank. 

\fphmg{} displays improved \relTotal{} in 3D compared to the 2D case. Specifically, 3D \fphmg{} results show faster time to solution even for the smallest problem types, achieving  1.5--3.2x speed-up across all problems relative to \fhmg{}.  
For lower-order discretization, \relTotal{} values are proportional to the \relSolve{} time as the solve phase dominates the run-time. For high discretization orders, run-time becomes dominated by setup cost, leading to larger improvements in \relTotal{} due to significantly faster setup times. 

The \phmg{} preconditioners, on the other hand, are all slower or match the total time to solution of \fhmg{}. This is attributed to much slower setup times due to bigger patches that contain a greater number of velocity DoFs (closure of vertex star as opposed to just vertex star) and the addition of pressure DoFs in the patch. 
\relSolve{} ratios for \phmg{} are 1.1--1.7 times faster than \fhmg{}, but they are  up to 2x slower than \fphmg{} for highest discretization orders. 
As in the setup phase, this issue is mainly attributed to relaxation cost on the $\ell=0$, Scott-Vogelius, grid. However, other kernels are marginally slower as well. 

It is important to point out relative improvement in~\relSolve{} for the \VphmgGradual{} solvers over \VphmgDirect{}. Despite this, we can see that the iteration counts for \VphmgGradual{} grow as well, albeit at a lower pace than \VphmgDirect{}. 

Alternative parameter combinations have been explored for \phmg{} solvers. In the case of \VphmgDirect{} and two refinement meshes problems, we observed a strong need for multiple \hmg{} cycles, suggesting a more robust low-order correction is necessary. 
The same applies for \VphmgGradual{} solvers, but in that case, additional \hmg{} cycles can be avoided using a more gradual $p$\kern1pt-coarsening schedule.
A useful heuristic we identified when selecting cycle parameters for the
Scott-Vogelius \phmg{} solvers is to reduce the amount of relaxation on the
Scott-Vogelius grid as it is significantly more costly than the relaxation on
the coarse grids. However, the communication cost associated with the coarse-grid
corrections is likely to become a bottleneck when the number of MPI tasks is
increased.

The number of unknowns and the average number of non-zero entries per degree of freedom (DoF) in the matrix operator for each problem shown in~\cref{fig:sv_ref2_ldc3d} are listed in~\cref{table:nnz_sv_3d}.
\begin{table}[htbp]
\centering
\begin{tabular}{c|c|rrrrrrrr}
\toprule
\multicolumn{2}{c|}{} & \multicolumn{8}{c}{\textbf{order} ($k$)} \\
\textbf{problem type} & \textbf{metric} & 3 & 4 & 5 & 6 & 7 & 8 & 9 & 10 \\
\midrule
\multirow{2}{*}{ldc3d} & DoFs (mil.) & 13 & 13 & 13 & 16 & 10 & 9 & 6 & 8 \\
                       & nnz/DoFs   & 132 & 220 & 337 & 487 & 673 & 899 & 1170 & 1491 \\
\bottomrule
\end{tabular}
\caption{Total number of degrees of freedom (DoFs, in millions) and the average number of non-zero entries per DoF in the Stokes operator discretized using the Scott-Vogelius element pair. The data is presented for different approximation orders ($k$) for 3D lid-driven cavity (ldc3D) problem with two mesh refinements (\cref{fig:sv_ref2_ldc3d}).}
\label{table:nnz_sv_3d}
\end{table}

\begin{remark}[Setup Costs]
It is important to acknowledge the high setup and storage costs involved in mixed-field patches. 
Exploring new research directions to reduce the costs associated with storage, setup, and application for patch relaxation techniques is paramount for advancing these methods and expanding their applicability. 
Some promising recent work in this area includes fast-diagonalization and patch compression techniques~\cite{harper2023compression,brubeck2021scalable,brubeck2022multigrid}.
Such advancements are key to future progress and a wider range of applications of these methods.
\end{remark}

\section{Lessons Learned}\label{sec:lessons_learned}

The examination of results for 2D and 3D problems discretized with Scott-Vogelius and Taylor-Hood elements, as discussed in~\cref{sec:th,sec:sv}, reveals several important insights regarding the performance and efficiency of multilevel solvers. These insights are essential for informed solver selection, understanding the trade-offs involved, and prioritizing future research efforts.

\subsection*{Taylor-Hood Discretizations}

\begin{itemize}
    \item In both \hmg{} and \phmg{} solvers, a significant portion of the runtime, often more than half, is dedicated to the setup phase. The fastest solvers, in terms of both setup and solution times, are those that can efficiently perform the setup phase. Given that both solver types use the same relaxation method at the finest level, which predominates setup cost, advancements in setup phase performance are primarily achieved through cheaper coarser level corrections.
    \item \phmg{} solvers improve efficiency at higher discretization orders through the use of sparser coarse grids and more cost-effective coarse-grid relaxation. This highlights the efficacy of polynomial coarsening strategies. However, as observed in \VphmgDirect{}, overly aggressive coarsening can deteriorate convergence, offsetting the benefits of cheaper coarse-grid corrections.
    \item The choice between direct and gradual $p$\kern1pt-coarsening strategies in
        \phmg{} solvers is important, with gradual coarsening providing better
        convergence rates but at a higher per-iteration cost and increased setup
        and memory overhead.
    \item Grid-transfer costs increase significantly with discretization order in \hmg{} solvers, an issue not present in \phmg{} solvers.
    \item For 3D problems, \phmg{} solvers show improved setup and solution phase timings relative to \hmg{} solvers. 
\end{itemize}

\subsection*{Scott-Vogelius Discretizations}

\begin{itemize}
    \item Similar to the findings for Taylor-Hood discretizations, \phmg{} solvers spend a considerable portion of the total runtime in the setup phase, emphasizing the need for optimization in this area.
    \item In 3D problems, both \afbf{} and \phmg{} solvers experience a significant increase in setup phase duration as a fraction of total computational time, particularly at higher discretization orders, which corresponds to rapid growth in matrix-stencil and thereby patch sizes.
    \item \fphmg{} solvers excel over both \fhmg{} and the \phmg{} solvers, especially at higher discretization levels, due to more efficient setup and solve phase timings. This underscores the benefit of integrating full-block factorization with \phmg{} on nested mesh hierarchies.
    \item \fphmg{} solvers outperform \fhmg{} across all discretization orders in 3D, benefiting from lower per-iteration costs stemming from more efficient coarse-grid relaxation and cheaper grid transfers.
    \item Achieving low iteration counts with \phmg{} solvers in 3D requires either increased V-cycle complexity (number of \hmg{} cycles, relaxation sweeps) or a more gradual $p$\kern1pt-coarsening strategy. The balance between these approaches and the cost of fine-level relaxation dominating runtime must be carefully considered.
\end{itemize}

\section{Conclusion}\label{sec:conclusion}

This paper introduced spatial ($h$) and approximation order ($p$)  multigrid (\phmg{}) methods for solving high-order discretized stationary Stokes systems employing Taylor-Hood and Scott-Vogelius elements.
A detailed comparative analysis showed that the \phmg{} method offers significant computational and memory efficiency improvements over traditional spatial-coarsening-only multigrid (\hmg) approaches, especially for problems discretized with Taylor-Hood elements. 
The effectiveness of this method is highlighted by its superior performance in reducing setup and solve times, particularly in scenarios involving higher discretization orders and unstructured meshes.

In the context of Scott-Vogelius discretizations, the comparison between monolithic \phmg{} and multilevel full-block-factorization (\afbf{}) preconditioners reveals that while monolithic \phmg{} achieves competitive solve phase timings, it is notably less efficient during the setup phase.
This efficiency gap is primarily due to the novel integration of \phmg{} within the \afbf{} preconditioner framework, which utilizes nested-mesh hierarchies. 
The \afbf{} preconditioner's advantage in performance largely stems from its lower setup costs for patch relaxation, attributed to involving a single unknown type. This contrasts with the monolithic \phmg{} method, which necessitates assembling larger mixed-field relaxation patches, making the setup phase more costly in comparison.

\section{Future Directions}\label{sec:future_directions}

The outcomes of this research highlight several promising directions for
advancing the capabilities of the \phmg{} framework. 
Key areas for future research include:
\begin{itemize}
	\item \textbf{Optimization of setup phase for \phmg{}}: Addressing the high setup costs associated with \phmg{} %
	presents a significant opportunity for improvement. Future work could explore more efficient algorithms for assembling mixed-field relaxation patches or alternative strategies that reduce memory requirements and the computational burden during the setup and solve phases~\cite{brubeck2021scalable,harper2023compression}
	\item \textbf{Optimization for Lower-Order Discretizations and Grids}: Enhancing \phmg{} for lower-order discretizations offers a significant opportunity to extend efficiency and applicability. Adopting alternative coarse-grid relaxation strategies, coupled with the potential integration of monolithic algebraic multigrid (AMG) methods as proposed in~\cite{voronin2023monolithic}, can result in notable performance improvements.
	\item \textbf{Extension to other coupled systems}: The promising results obtained for stationary Stokes systems suggest that the \phmg{} framework could be adapted and applied to other complex systems, including those involving non-linear equations or additional physics.
\end{itemize}

\section*{Acknowledgements}
The work of SPM was partially supported by an NSERC Discovery Grant. RT was
supported by the U.S.~Department of Energy, Office of Science, Office of
Advanced Scientific Computing Research, Applied Mathematics program. Sandia
National Laboratories is a multi-mission laboratory managed and operated by
National Technology \& Engineering Solutions of Sandia, LLC (NTESS), a wholly owned subsidiary of Honeywell International Inc., for the U.S. Department of Energy’s National Nuclear Security Administration (DOE/NNSA) under contract DE-NA0003525. This written work is authored by an employee of NTESS. The employee, not NTESS, owns the right, title and interest in and to the written work and is responsible for its contents. Any subjective views or opinions that might be expressed in the written work do not necessarily represent the views of the U.S. Government. The publisher acknowledges that the U.S. Government retains a non-exclusive, paid-up, irrevocable, world-wide license to publish or reproduce the published form of this written work or allow others to do so, for U.S. Government purposes. The DOE will provide public access to results of federally sponsored research in accordance with the DOE Public Access Plan.

This work used the Delta system at the National Center for Supercomputing Applications through allocation CIS230361 from the Advanced Cyberinfrastructure Coordination Ecosystem: Services \& Support (ACCESS) program~\cite{boerner2023access}, which is supported by National Science Foundation grants \#2138259, \#2138286, \#2138307, \#2137603, and \#2138296.

\bibliographystyle{siamplain}
\bibliography{references}
\end{document}